\documentclass[a4paper,peperxidth=110mm,paperheigth=110mm,10pt,leqno]{article}
\usepackage{amsmath}
\usepackage{amsfonts}
\usepackage{hyperref}
\hypersetup{
    colorlinks=true,                         
    linkcolor=red, % Couleur des liens internes
    citecolor=blue, % Couleur des numéros de la biblio dans le corps
    urlcolor=red  } % Couleur des url
\usepackage{graphicx}
\usepackage{latexsym}
\usepackage{fancybox}
\usepackage[utf8]{inputenc}
\usepackage{amssymb}
\usepackage{fancyhdr}
\usepackage{amsthm}
\usepackage{tikz}
\usepackage[a4paper,tmargin=1.5cm,bmargin=1.5cm,
rmargin=1cm,lmargin=1cm]{geometry}
\newtheorem{theorem}{Theorem}[section]

\newtheorem{lemma}{Lemma}[section]
\newtheorem{rem}{Remark}[section]

\newtheorem{proposition}{Proposition}[section]

\newenvironment{demo}{\noindent{Proof of \noindent}}{\hfill  $\Box$\vskip 5mm}

\numberwithin{equation}{section}

\begin{document}

\title{Null controllability of a nonlinear age and two-sex population dynamics structured model}
\author{Amidou Traore \footnote{Laboratoire LAMI, Université Joseph Ki ZERBO, 01 BP 7021 Ouaga 01, Burkina Faso (\texttt{ amidoutraore70@yahoo.fr})}
\and
Okana S. Sougué \footnote{Laboratoire LAMI, Université Joseph Ki ZERBO, 01 BP 7021 Ouaga 01, Burkina Faso (\texttt{ sougueok@gmail.com})}
\and
 Yacouba Simporé
\footnote{Laboratoire LAMI, Université Joseph Ki ZERBO, 01 BP 7021 Ouaga 01, Burkina Faso. DeustoTech, Fundación Deusto Avda. Universidades,
24, 48007, Bilbao, Basque Country, Spain  (\texttt{simplesaint@gmail.com})}
 \and
 Oumar Traore\footnote{Département de Mathématiques de la Décision, Université Ouaga 2, 12 BP 417 Ouaga 12, Laboratoire LAMI, Université Joseph Ki ZERBO, 01 BP 7021 Ouaga 01,
Burkina Faso (\texttt{otraore.@univ-ouaga2.bf})}}
\date{}
\maketitle

\begin{abstract}
This paper is devoted to study the null controllability properties of a nonlinear age and two-sex population dynamics structured model without spatial structure. Here, the nonlinearity and the couplage are at birth level.

\noindent In this work we consider two cases of null controllability problem. 

\noindent  The first problem is related to the extinction of male and female subpopulation density. 

\noindent The second case concerns the null controllability of male or female subpopulation individuals. In both cases, if $A$ is the maximal age expectancy, a time interval of duration $A$ after the extinction of males or females, one must get the total extinction of the population.  

\noindent Our method uses first an observability inequality related to the adjoint of an auxiliary system, a null controllability of the linear auxiliary system and after the Kakutani's fixed point theorem.  
\end{abstract}

\section{Introduction and main results}
In this paper, we study the null controllability of an infinite dimensional nonlinear coupled system describing the dynamics of two-sex structured population without spatial position. Let $(m,f)$ be the solution of the following system :
\begin{equation}\label{EP}
\left\lbrace\begin{array}{ll}
m_t+m_a+\mu_m m  =\chi_\Xi v_m &\hbox{in}\;Q, \\
\\ 
f_t+f_a+\mu_f f  =\chi_{\Xi'} v_f &\hbox{in}\;Q, \\
\\ 
m(a,0)=m_0(a)\quad f(a,0)=f_0(a) & \hbox{in}\;Q_A, \\
\\ 
m(0,t)=(1-\gamma)\displaystyle\int_0^A\beta(a,M)f(a,t)da & \hbox{in}\;Q_T, \\
\\ 
f(0,t)=\gamma\displaystyle\int_0^A\beta(a,M)f(a,t)da & \hbox{in}\;Q_T, \\ 
\\
M=\displaystyle\int_0^A\lambda(a)m(a,t)da & \hbox{in}\;Q_T,
\end{array}\right. 
\end{equation} 
\noindent where $T$ is a positive number, $Q=(0,A)\times(0,T),\;\Theta=(0,a_2)\times(0,T),\;\Xi=(a_1,a_2)\times(0,T)$ and $\Xi'=(b_1,b_2)\times(0,T).$ Here $0\leq a_1<a_2\leq A,\;0\leq b_1<b_2\leq A,\;Q_A=(0,A)$ and $Q_T=(0,T).$

\noindent We denote the density of males and females of age $a$ at time $t$ respectively by $m(a,t)$ and $f(a,t).$ Moreover, $\mu_m$ and $\mu_f$ denote respectively the natural mortality rate of males and females. The control functions are $v_m$ and $v_f$ and depend on $a$ and $t.$ In addition $\chi_\Xi$ and $\chi_{\Xi'}$ are the characteristic functions of the support of the control $v_m$ and $v_f$ respectively.

\noindent We have denoted by $\beta$ the positive function describing the fertility rate that depends on $a$ and also on 
\[M=\displaystyle\int_0^A\lambda(a)m(a,t)da,\]
\noindent where $\lambda$ is the fertility function of the male individuals. Thus the densities of newborn male and female individuals at time $t$ are given respectively by $m(0,t)=(1-\gamma)N(t)$ and $f(0,t)=\gamma N(t)$ where
\[N(t)=\displaystyle\int_0^A\beta(a,M)f(a,t)da.\]
\noindent We assume that the fertility rate $\beta,\;\lambda$ and the mortality rate $\mu_f,\;\mu_m$ satisfy the demographic properties :
\begin{align*}
(H_1) &\left\lbrace\begin{array}{l}
\mu_m(a)\geq 0,\quad \mu_f(a)\geq 0\;\hbox{a.e}\;a\in(0,A) \\ 
\\
\mu_m\in L^1_{loc}(0,A),\quad \mu_f\in L^1_{loc}(0,A)\\ 
\\
\displaystyle\int_0^A\mu_m(a)da=+\infty,\quad \displaystyle\int_0^A\mu_f(a)da=+\infty
\end{array} \right.\\
\\
(H_2) & \left\lbrace\begin{array}{l}
\beta(a,p)\in C([0,A]\times\mathbb{R}) \\
\\ 
\beta(a,p)\geq 0\;\hbox{for every}\;(a,p)\in[0,A]\times\mathbb{R}.
\end{array} \right.
\end{align*}
\noindent We further assume that the birth function $\beta$ and the fertility function $\lambda$ verify the following hypothesis:
\begin{align*}
(H_3)& \left\lbrace\begin{array}{l}
\hbox{there exists}\;b\in(0,A)\;\hbox{such that}\;\beta(a,p)=0,\;\forall(a,p)\in(0,b)\times\mathbb{R}, \\ 
\\
\hbox{there exists}\;\vert\vert\beta\vert\vert_\infty>0\;\hbox{such that}\;0\leq \beta\leq \vert\vert\beta\vert\vert_\infty,\;\forall(a,p)\in(0,b)\times\mathbb{R},
\\
 \\ 
\beta(a,0)=0,\;\forall\;a\in (0,A).
\end{array} \right.\\
\\
(H_4) &\left\lbrace\begin{array}{l}
\lambda\in C^1([0,A]) \\
\\ 
\lambda(a)\geq 0\;\hbox{for every}\;a\in[0,A],\\
\\
\lambda\mu_m\in L^1(0,A).
\end{array} \right.
\end{align*}
\noindent The assumption $\beta(a,0)=0$ for $a\in(0, A)$ means that, the birth rate is zero if there are no fertile male individuals.

\noindent We can now state the main results. If $(a_1,a_2)\subset(b_1,b_2),$ we have the following theorem:
\begin{theorem}\label{theo1}
Let us assume that the assumptions $(H_1)-(H_4)$ hold true. If $a_1<b,$ for every time 

\noindent $T>a_1+A-a_2$ and for every $(m_0,f_0)\in \big(L^2(Q_A)\big)^2,$ there exists $(v_m,v_f)\in L^2(\Xi)\times L^2(\Xi')$ such that the associated solution $(m,f)$ of system $(\ref{EP})$ verifies:
\begin{equation}\label{ezero1}
m(a,T)=f(a,T)=0\quad\hbox{a.e}\;a\in(0,A).
\end{equation}
\end{theorem}

\begin{theorem}\label{theo2}
Let us assume that the assumptions $(H_1)-(H_4)$ hold true. We have :
\begin{enumerate}
\item[(1)] let $v_f=0.$ For any $\varrho>0,$ for every time $T>A-a_2$ and for every $(m_0,f_0)\in \big(L^2(Q_A)\big)^2,$ there exists a control $v_m\in L^2(\Theta)$ such that the associated solution $(m,f)$ of system $(\ref{EP})$ verifies:
\begin{equation}\label{ezero2}
m(a,T)=0\quad\hbox{a.e}\;a\in(\varrho,A)
\end{equation}
\noindent where $\Theta=(0,a_2)\times(0,T).$
\item[(2)] let $v_m=0.$ For every time $T>a_1+A-a_2$ and for every $(m_0,f_0)\in \big(L^2(Q_A)\big)^2,$ there exists a control $v_f\in L^2(\Xi)$ such that the associated solution $(m,f)$ of system $(\ref{EP})$ verifies:
\begin{equation}\label{ezero2}
f(a,T)=0\quad\hbox{a.e}\;a\in(0,A).
\end{equation}
\end{enumerate}
\end{theorem}

\begin{rem}
The first condition of $(H_3)$ is not necessary for the Theorem~\ref{theo2}-(1).
\end{rem}
In practice this study takes place in the fight against malaria. Malaria is a serious disease (in $2017$ there were $219$ million cases in the world \cite{WHO}) and our work takes its importance in the strategy to fight against it.

\noindent In fact Malaria is a vector-borne disease transmitted by an infective female anopheles mosquito. A malaria control strategy
in Brazil or Burkina Faso consists of releasing genetically modified male mosquitoes (precisely sterile males) in the nature.
This can reduce the reproduction of mosquito since females mate only once in their life cycle.

\noindent In the theoretical framework, very few authors have studied control problems of two-sex structured population dynamics
model.

\noindent The control problems of coupled systems of population dynamics models take an intense interest and are widely investigated in many papers. Among them we have \cite{BASA2}, \cite{MIJP}, \cite{CZMWPZ} and the references therein. In fact, in \cite{BASA2} the authors studied a coupled
reaction-diffusion equations describing interaction between a prey population and a predator one. The goal of the above work is to look for a suitable control supported on a small spatial subdomain which guarantees the stabilization of the predator
population to zero. In \cite{CZMWPZ}, the objective was different. More precisely, the authors consider an age-dependent prey predator system and they prove the existence and uniqueness of an optimal control (called also "optimal effort") which gives the maximal harvest via the study of the optimal harvesting problem associated to their coupled model.
In \cite{YHBA} He and Ainseba study the null controllability of a butterfly population by acting on eggs, larvas and female moths
in a small age interval. 

\noindent In \cite{MIJP}, the authors analyze the growth of a two-sex population with a fixed age-specific sex ratio without diffusion. The
model is intended to give an insight into the dynamics of a population where the mating process takes place at random
choice and the proportion between females and males is not influenced by environmental or social factors, but only depends
on a differential mortality or on a possible transition from one sex to the other (e.g. in sequential hermaphrodite species). Simpore and Traore study in \cite{YSOT} the null controllability of a nonlinear age, space and two-sex structured population dynamics model. They first study an approximate null controllability result for an auxiliary cascade system and prove the null controllability of the nonlinear system by means of Schauder’s fixed point theorem. 

\noindent Unlike the model treated in \cite{YSOT}, we consider a nonlinear cascade system with two different fertility rates and without space variable. The fertility rate of the male $\lambda$ and the fertility rate $\beta$ of the female depend on the total population of the fertile males.

\noindent We use the technique of \cite{YSOT} and \cite{DM} combining final-state observability estimates with the use of characteristics to establish the observability inequalities necessary for the null controllability property of the auxiliary systems. Roughly, in our method we first study the null controllability result for an auxiliary cascade system. Afterwards, we prove the null controllability result for the system $(\ref{EP})$ by means of Kakutani’s fixed point theorem.

\noindent The remainder of this paper is as follows: in Section $2$ we study the existence and uniqueness of a positive solution for the model. Section $3$ is devoted to the proof of the Theorem~\ref{theo1} and the Theorem~\ref{theo2}.

\section{Well posedness result}\label{section1}

 In this section, we study the existence of positive solution of the model. For this, we assume that the so-called demographic conditions $(H_1),\;(H_2),\;(H_3)$ and $(H_4)$ are verified. Moreover, here, we suppose that 
\begin{align*}
 &(H_5)\left\lbrace\begin{array}{l}
 \beta(a,p)=\beta_1(a)\beta_2(p)\;\hbox{for all}\;(a,p)\in(0,A)\times\mathbb{R}, \\
 \\ 
\hbox{there exists}\;C>0\;\hbox{such that}\;\vert\beta_2(p)-\beta_2(q)\vert\leq  C\vert p-q\vert\;\hbox{for all}\; p,q\in\mathbb{R},\\
\\ 
 \beta_1\mu_f\in L^1(0,A)
 \end{array}\right. 
\end{align*}  
\noindent holds true.

\noindent Thus, we have the following result.
\begin{theorem}\label{theo3}
$\;$

\noindent Assume that $(H_1)-(H_5)$ hold. For every $(m_0,f_0)\in(L^2(0,A))^2$ and $(v_m,v_f)\in(L^2(Q))^2,$ the system $(\ref{EP})$ admits a unique solution $(m,f)\in(L^2((0,A)\times(0,T)))^2$ and we have the following estimations:
\begin{equation}\label{ees1}
\left\lbrace\begin{array}{l}
\vert\vert m\vert\vert_{L^2((0,A)\times(0,T))}\leq K\big(\vert\vert f_0\vert\vert_{L^2(0,T)}+\vert\vert m_0\vert\vert_{L^2(0,T)}+\vert\vert v_m\vert\vert_{L^2(Q)}+\vert\vert v_f\vert\vert_{L^2(Q)}\big), \\
\\ 
\vert\vert f\vert\vert_{L^2((0,A)\times(0,T))}\leq C\big(\vert\vert m_0\vert\vert_{L^2(0,T)}+\vert\vert v_f\vert\vert_{L^2(Q)}\big)
\end{array}\right. 
\end{equation}
\noindent where $K$ and $C$ are positive constants.

\noindent Moreover, suppose that 
\[ 
m_0,f_0\geq 0\;\hbox{a.e}\;(0,A)\;\hbox{and}\;v_m,v_f\geq 0\;\hbox{a.e}\;Q;
 \]
\noindent then $(m,f)$ is also positive. 
\end{theorem} 
\begin{demo}Theorem~\ref{theo3}

\noindent Let $p$ be fixed in $L^2(0,T)$, $h$ and $h'$ be fixed in $L^2(Q)$ and consider the following system
\begin{equation}\label{Ok}
\left\lbrace\begin{array}{ll}
m_t+m_a+\mu_m m  =h &\hbox{in}\;Q, \\
\\ 
f_t+f_a+\mu_f f  =h' &\hbox{in}\;Q, \\
\\ 
m(a,0)=m_0(a)\quad f(a,0)=f_0(a) & \hbox{in}\;Q_A, \\
\\ 
m(0,t)=(1-\gamma)\displaystyle\int_0^A\beta\left(a,\int_0^A \lambda (a)p(a,t)\right)f(a,t)da & \hbox{in}\;Q_T, \\
\\ 
f(0,t)=\gamma\displaystyle\int_0^A\beta\left(a,\int_0^A \lambda(a)p(a,t)\right)f(a,t)da & \hbox{in}\;Q_T, \\ 
\end{array}\right. 
\end{equation} 
\noindent For every $f_0\in L^2(0,A)$ and $h'\in L^2(Q)$, the following system
\begin{equation}\label{Ok}
\left\lbrace\begin{array}{ll}
f_t+f_a+\mu_f f  =h' &\hbox{in}\;Q, \\
\\ 
 f(a,0)=f_0(a) & \hbox{in}\;Q_A, \\
\\
f(0,t)=\gamma\displaystyle\int_0^A\beta\left(a,\int_0^A \lambda(a)p(a,t)\right)f(a,t)da & \hbox{in}\;Q_T, \\ 
\end{array}\right. 
\end{equation} 
\noindent admits a unique positive solution in $L^2(Q),$ see \cite{SA} and one has 
 \begin{equation}\label{Ineq}
 \left\|f\right\|^2_{L^2(Q)}\leqslant C\left( \left\|f_0\right\|^2_{L^2(0,A)}+\left\|h'\right\|^2_{L^2(Q)}\right),
 \end{equation}
\noindent where $C$ is a positive constant and independent of $p$ because $\beta\in L^\infty((0,T)\times(0,A)).$  
 
\noindent Now, $f$ and $h'$ are being known, the system
 \begin{equation}\label{AA}
\left\lbrace\begin{array}{ll}
m_t+m_a+\mu_m m  =h &\hbox{in}\;Q, \\
\\ 
 m(a,0)=m_0(a) & \hbox{in}\;Q_A, \\
\\
m(0,t)=(1-\gamma)\displaystyle\int_0^A\beta\left(a,\int_0^A \lambda(a)p(a,t)\right)f(a,t)da & \hbox{in}\;Q_T, \\ 
\end{array}\right.
\end{equation} 
\noindent admits a unique positive system in $L^2(Q)$ and we have the following estimation 
\[ \left\|m\right\|^2_{L^2(Q)}\leqslant K\left( \left\|f_0\right\|^2_{L^2(0,A)}+\left\|m_0\right\|^2_{L^2(0,A)}+
\left\| h\right\|^2_{L^2(Q)}+\left\|h'\right\|^2_{L^2(Q)}\right),\]
\noindent where $K$ is a positive constant and  independent of $p$ because $\beta\in L^\infty((0,T)\times(0,A)).$  

\noindent Let define $\Phi \;:\, L^2_+(Q)\longrightarrow L^2_+(Q)$, $\Phi p=m(p)$ where $m(p)$ is the unique solution of the system $(\ref{AA})$. 

\noindent For any $p\;,\; q \in L^2_+(Q)$, we set $$B_1(a,t)=\int_0^{A}\lambda (a)p(t,a)da\quad   \text{and}\quad B_2(a,t)=\int_0^{A}\lambda(a)q(t,a)\quad \text{a.e.}\quad t\in (0,A)\times(0,T),$$  and $w=(m(p)-m(q))e^{-\gamma_0t}$ where $\gamma_0$ is a positive parameter that will be choosed later; $w$ is solution of
\begin{equation}\label{kk}
\left\lbrace\begin{array}{ll} 
w_t+w_a+(\gamma_0+\mu_m) w  =0 &\hbox{in}\;Q, \\
\\ 
 w(a,0)=0 & \hbox{in}\;Q_A, \\
\\
w(0,t)=(1-\gamma)e^{-\gamma_0t} \times \\  \makebox[1cm][1cm]\,\displaystyle\int_0^A \left[\beta_2(B_1)-\beta_2(B_2)\right]\beta_1(a)f(p)+(f(p)-f(q)) \beta_2(B_2)\beta_1(a)da & \hbox{in}\;Q_T, \\ 
\end{array}\right.
\end{equation} 
\noindent Multipliying $(\ref{kk})$ bye $w$ and integrating over $(0,A)\times (0,t)$, and using Young's inegality we get 
\begin{align*}
\frac{1}{2}\left\| w(t) \right \|^2_{L^2(0,A)}+\int_0^t \int_0^A (\gamma_0+\mu_m)w^2(s,a)dads &\leqslant \int_0^t \left(\int_0^A\left| \beta_2(B_1)-\beta_2(B_2)\right| \beta_1(a)f(p)da \right)^2ds\\
&\quad +\int_0^t \left(\int_0^A\left(f(p)-f(q)\right)\beta_2(B_2)\beta_1(a)da\right)^2ds \\
&\leqslant C^2 \left\| \lambda\right\|^2_{\infty} \int_0^t \left(\left| \int_0^A p(s,a)da-\int_0^A q(s,a)da \right| \right)^2\left(\int_0^A \beta_1(a)f(p(s))da\right)^2ds \\ 
&\quad +\int_0^t \left(\int_0^A\left(f(p)-f(q)\right)\beta_2(B_2)\beta_1(a)da\right)^2ds \\
&\leqslant  C^2 A \left\| \lambda\right\|^2_{\infty} \int_0^t \int_0^A  \left|  p(s)- q(s) \right| ^2\left(\int_0^A \beta_1(a)f(p(s))da\right)^2ds\\
&\quad +\left\| \beta_1\right\|^2_{\infty} \left\|\beta_2\right\|^2_{\infty}A\int_0^t\int_0^A \left|f(p)-f(q)\right |^2ds. 
\end{align*}
\noindent Hence for every $\gamma_0>0$, there is a constant $C= \max \left\lbrace 2C^2 A \left\| \lambda\right\|^2_{\infty};2\left\| \beta_1\right\|^2_{\infty} \left\|\beta_2\right\|^2_{\infty}A\right\rbrace$ such that  \begin{equation}\label{Ine1}
\left\| w(t) \right \|^2_{L^2(0,A)}\leqslant C\left( \int_0^t \int_0^A  \left|  p(s)- q(s) \right| ^2\left(\int_0^A \beta_1(a)f(p(s))da\right)^2ds +\int_0^t\int_0^A \left|f(p(s))-f(q(s))\right |^2ds \right) 
\end{equation}
\noindent Now set $F =(f(p)-f(q))e^{-\delta t}$ where $\delta $ is a positive parameter that will be choosed later. Then, $F $ solves  the following auxiliary system
\begin{equation}\label{kkk}
\left\lbrace\begin{array}{ll} 
F_t+F_a+(\delta+\mu_f) F  =0 &\hbox{in}\;Q, \\
\\ 
 F(a,0)=0 & \hbox{in}\;Q_A, \\
\\
w(0,t)=\gamma \int_0^A e^{-\gamma_0t}\left[\beta_2(H_1)-\beta_2(H_2)\right]\beta_1(a)f(p)+F(a,t) \beta_2(H_2)\beta_1(a)da & \hbox{in}\;Q_T. \\ 
\end{array}\right.
\end{equation}
\noindent Similarly as above, we have 
\[\delta \int_0^t\int_0^A F(a,s)^2 ds \leqslant C\left( \int_0^t \int_0^A  \left|  p(s)- q(s) \right| ^2\left(\int_0^A \beta_1(a)f(p(s))da\right)^2ds +\int_0^t\int_0^A \left|F(a,s)\right |^2ds \right).\]
\noindent Hence, there is a positive constant $C'$ such that 
\begin{equation}\label{Ine2}
\int_0^t\int_0^A F(a,s)^2 ds \leqslant C' \int_0^t \int_0^A  \left|  p(s)- q(s) \right| ^2\left(\int_0^A \beta_1(a)f(p(s))da\right)^2ds .
\end{equation}
\noindent Setting $Y(t)=\displaystyle\int_0^A\beta_1(a)f(p)da\; \text{a.e}\; \text{in}\; (0, T)$, $Y$ solves the following system 
\begin{equation}\label{kkkk}
\left\lbrace\begin{array}{ll} 
Y_t  =\displaystyle\int_0^A\beta_1(a)h'(a,t)da+\displaystyle\int_0^A\beta_1'(a)f(a,t)da-\displaystyle\int_0^A\mu_f(a)\beta_1(a)f(a,t)da &\hbox{in}\;(0,T),\\
\\ 
 Y(0)= \displaystyle\int_0^A \beta_1(a)f_0(a)da  \\
\end{array}\right.
\end{equation}
\noindent Multiplying $(\ref{kkkk}) $ by $Y,$ integrating over $(0,t)$ and using Young's inegality we get  
\begin{align*}
 Y^2(t) &\leqslant Y^2(0) +\int_0^tY^2(s)ds +\int_0^t\left(\int_0^A\beta_1(a)h'(a,s)da+\int_0^A \beta_1'(a)f(a,s)da+\int_0^A \mu_f(a)\beta_1(a)f(a,s)da\right)^2ds \\
 &\leqslant Y^2(0)+\int_0^t Y^2(s) ds+3 \int_0^t\left(\int_0^A \beta_1(a)h'(a,s)da\right)^2ds +3\int_0^t \left(\int_0^A \beta_1'(a)f(a,s)da\right)^2ds\\
 &\quad +3\int_0^t \left(\int_0^A \beta_1(a)\mu_f(a)f(a,s)da\right)^2ds. 
\end{align*}
\noindent So,
\begin{align}\label{Ok2}
Y^2(t)&\leqslant \left(\int_0^A \beta_1(a)f_0(a)da\right)^2+\int_0^T Y^2(t) dt +3 \int_0^T\left(\int_0^A \beta_1(a)h'(a,t)da\right)^2dt +3\int_0^T \left(\int_0^A \beta_1'(a)f(a,t)da\right)^2dt\\
 &\quad +3\int_0^T \left(\int_0^A \beta_1(a)\mu_f(a)f(a,t)da\right)^2dt. \nonumber
\end{align}
\noindent Let us set $\tilde{f}=e^{-\lambda_0 t}f.$ Then, from $(\ref{Ok})\;\tilde{f}$ satisfies the following system
\begin{equation}\label{Ok1}
\left\lbrace\begin{array}{ll}
\tilde{f}_t+\tilde{f}_a+(\lambda_0+\mu_f) \tilde{f}  =e^{-\lambda_0 t}h' &\hbox{in}\;Q, \\
\\ 
 \tilde{f}(a,0)=f_0(a) & \hbox{in}\;Q_A, \\
\\
\tilde{f}(0,t)=\gamma\displaystyle\int_0^A\beta\left(a,\int_0^A \lambda(a)p(a,t)\right)\tilde{f}(a,t)da & \hbox{in}\;Q_T, \\ 
\end{array}\right. 
\end{equation}
\noindent Multiplying the first equation of $(\ref{Ok1})$ by $\tilde{f},$ integrating on $Q$ and using the inequality of Young we get
\begin{align*}
& \int_0^T \int_0^A(\lambda_0+\mu_f(a))\tilde{f}^2(a,t)dadt\leq \dfrac{1}{2}\vert\vert f_0\vert\vert^2_{L^2(0,A)}+\dfrac{1}{2}\vert\vert h'\vert\vert^2_{L^2(Q)}+\dfrac{1}{2}\vert\vert \tilde{f}\vert\vert^2_{L^2(Q)}+\dfrac{1}{2}\int_0^T\tilde{f}^2(0,t)dt.
\end{align*}
\noindent Using the inequality of Cauchy Schwartz and choosing $\lambda_0=\dfrac{3}{2}+\vert\vert \beta\vert\vert^2_{\infty},$ we obtain
\begin{align*}
& \int_0^T \int_0^A \mu_f(a)\tilde{f}^2(a,t)dadt\leq \dfrac{1}{2}\bigg(\vert\vert f_0\vert\vert^2_{L^2(0,A)}+\vert\vert h'\vert\vert^2_{L^2(Q)}\bigg).
\end{align*}
\noindent So,
\begin{align}\label{Ok3}
& \int_0^T \int_0^A \mu_f(a)f^2(a,t)dadt\leq \dfrac{e^{(3+2\vert\vert \beta\vert\vert^2_{\infty} )T}}{2}\bigg(\vert\vert f_0\vert\vert^2_{L^2(0,A)}+\vert\vert h'\vert\vert^2_{L^2(Q)}\bigg).
\end{align}
\noindent Using $(\ref{Ok2}),\;(\ref{Ok3})$ and against Young's inegality we have 
\begin{align*}
& Y^2(t) \leqslant \left\|\beta_1\right\|_{\infty}^2A\left\| f_0\right\|^2_{L^2(0,A)}+\left\|\beta_1\right\|_{\infty}^2A\left\| f\right\|^2_{L^2(Q)} +3 \left\|\beta_1\right\|_{\infty}^2A\left\| h'\right\|^2_{L^2(Q)}\\
& \qquad\qquad+3\left\|\beta'_1\right\|_{\infty}^2A\left\| f\right\|^2_{L^2(Q)} +3C\vert\vert\beta_1\vert\vert_\infty\left\|\beta_1 \mu_f\right\|_{L^1(0,A)}\vert\vert f_0\vert\vert^2_{L^2(0,A)}\\
&\qquad\qquad +3C\vert\vert\beta_1\vert\vert_\infty\left\|\beta_1 \mu_f\right\|_{L^1(0,A)}\vert\vert h'\vert\vert^2_{L^2(Q)}.
\end{align*}
\noindent From $(\ref{Ineq})$, we have just proved the existence of a positive constant $C$ such that , \begin{equation}\label{ineqq}
Y^2(t)\leqslant C\left(\left\|f_0\right\|^2_{L^2(0,A)}+\left\|h'\right\|^2_{L^2(0,A)}\right)
\end{equation} 
\noindent The estimate $(\ref{ineqq})$ means also that $Y\in L^{\infty}(0,T)$.

\noindent Combining $(\ref{Ine1})$, $(\ref{Ine2})$ and $(\ref{ineqq})$, we get the following estimate 
\begin{equation}\label{Ine3}
\left\| (\Phi p -\Phi q)(t) \right \|^2_{L^2(0,A)}\leqslant \sigma \int_0^t  \left\|  p(s)- q(s) \right\| ^2_{L^2(0,A)}ds,
\end{equation} 
\noindent where $\sigma $ is a positive constant.

\noindent Let us define the metric $d$ on $L^2_+(Q)$ by setting
\[d(h_1,h_2)=\left( \int_{0}^{T}\| (h_1-h_2)(t)\| ^2 _{L^2((0,A)}\exp\{-2\sigma t\}dt\right)^{\frac{1}{2}}, \quad \text{for} \; h_1, h_2 \in L^2_+(Q).\]
\noindent We have
\[d(\Phi p,\Phi q)^2 =\int_{0}^{T}\| (\Phi p-\Phi q)(t)\| ^2 _{L^2((0,A)}\exp\{-2\sigma t\}dt \leqslant \sigma\int_{0}^{T}\exp\{-2\sigma t\}\int_0^t \| ( p-q)(s)\| ^2 _{L^2((0,A)} ds dt\]
\noindent Using the Fubbini theorem, we conclude that
\begin{align*}
d(\Phi p,\Phi q)^2 =\int_{0}^{T}\| (\Phi p-\Phi q)(t)\| ^2 _{L^2((0,A)}\exp\{-2\sigma t\}dt &\leqslant \int_{0}^{T} \| ( p-q)(s)\| ^2 _{L^2((0,A)}\times \int_s^T \sigma e^{-2\sigma t}dtds  \\
 &\leqslant \frac{1}{2}d(p,q)^2.
\end{align*}
\noindent Then, $\Phi$ is a contraction on the complete metric space $L^2_+(Q)$  into itself. Using Banach's fixed point theorem,  we conclude the existence of a unique fixed point $m$. Moreover, $m$ is nonnegative. Hence, the unique couple $(m,f)$ is the unique solution to our problem $(\ref{EP})$.
\end{demo}

\section{Null controllability results}
\noindent For the sequel, the hypothesis $(H_5)$ is not necessary. 
\subsection{Null controllability of an auxiliary coupled system}

\noindent This section is devoted to the study of an auxiliary system obtained from the system $(\ref{EP})$.

\noindent Let $p$ be a $L^2(Q_T)$ function, we define the auxiliary system given by: 
\begin{equation}\label{EPauxiliary}
\left\lbrace\begin{array}{ll}
m_t+m_a+\mu_m m  =\chi_\Xi v_m &\hbox{in}\;Q, \\
\\ 
f_t+f_a+\mu_f f  =\chi_{\Xi'} v_f &\hbox{in}\;Q, \\
\\ 
m(a,0)=m_0(a)\quad f(a,0)=f_0(a) & \hbox{in}\;Q_A, \\
\\ 
m(0,t)=(1-\gamma)\displaystyle\int_0^A\beta(a,p)f(a,t)da & \hbox{in}\;Q_T, \\
\\ 
f(0,t)=\gamma\displaystyle\int_0^A\beta(a,p)f(a,t)da & \hbox{in}\;Q_T.
\end{array}\right. 
\end{equation} 
\noindent Let $p$ be fixed in $L^2(0,T),$ for $(m_0,f_0)\in\big(L^2(Q_A)\big)^2$ and $(v_m,v_f)\in L^2(\Xi)\times L^2(\Xi')$ the system $(\ref{EPauxiliary})$ admits a unique solution $(m,f)\in\big(L^2(Q)\big)^2,$ see Section~\ref{section1}.

\noindent The system above is null approximately controllable. Indeed we have the following result:
\begin{theorem}\label{theocontrolauxiliary}
Let us assume that assumptions $(H_1)-(H_2)$ hold. For every time $T>a_1+A-a_2,$ for every $\kappa,\theta>0$ and for every $(m_0,f_0)\in \big(L^2(Q_A)\big)^2,$ there exists a control $(v_\kappa,v_\theta)$ such that the solutions $m$ and $f$ of the system~$(\ref{EPauxiliary})$ verify
\begin{equation}\label{ebornitudesolution}
\vert\vert m(.,T)\vert\vert_{L^2(0,A)}\leq \kappa\quad\hbox{and}\quad \vert\vert f(.,T)\vert\vert_{L^2(0,A)}\leq \theta.
\end{equation}
\end{theorem}
\noindent The adjoint system of $(\ref{EPauxiliary})$ is given by:
\begin{equation}\label{SAdjoint}
\left\lbrace\begin{array}{ll}
-n_t-n_a+\mu_m n  =0 &\hbox{in}\;Q, \\
\\ 
-l_t-l_a+\mu_f l  =(1-\gamma)\beta(a,p)n(0,t)+\gamma\beta(a,p)l(0,t) &\hbox{in}\;Q, \\
\\ 
n(a,T)=n_T(a),\quad l(a,T)=l_T(a) & \hbox{in}\;Q_A, \\
\\ 
n(A,t)=0,\quad l(A,t)=0& \hbox{in}\;Q_T.
\end{array}\right. 
\end{equation} 
\noindent The main idea in this part is to establish an observability inequality of $(\ref{SAdjoint})$ that will allow us to prove the approximate null controllability of $(\ref{EPauxiliary}).$

\noindent The basic idea for establishing this inequality is the estimation of non-local terms. For that, let us start by formulating a representation of the solution of the adjoint system by the method of characteristic and semi-group.

\noindent For every $(n_T,l_T)\in (L^2(Q_A))^2,$ under the assumptions $(H_1)$ and $(H_2)$, the coupled system $(\ref{SAdjoint})$ admits a unique solution $(n,l).$ Moreover integrating along the characteristic lines, the solution $(n,l)$ of $(\ref{SAdjoint})$ is as follows:
\begin{equation}\label{esol1}
n(t)=\left\lbrace\begin{array}{l}
\dfrac{\pi_1(a+T-t)}{\pi_1(a)}n_T(a+T-t)\quad\hbox{if}\quad T-t\leq A-a, \\
\\ 
0 \quad\hbox{if}\quad A-a< T-t
\end{array} \right.
\end{equation}
\noindent and 
\begin{equation}\label{esol2}
l(t)=\left\lbrace\begin{array}{l}
\dfrac{\pi_2(a+T-t)}{\pi_2(a)}l_T(a+T-t)\\
+\displaystyle\int_t^T\frac{\pi_2(a+s-t)}{\pi_2(a)}\beta(a+s-t,p(s))\big((1-\gamma)n(0,s)+\gamma l(0,s)\big)ds\quad\hbox{if}\quad T-t\leq A-a, \\
\\ 
\displaystyle\int_t^{t+A-a}\frac{\pi_2(a+s-t)}{\pi_2(a)}\beta(a+s-t,p(s))\big((1-\gamma)n(0,s)+\gamma l(0,s)\big)ds \quad\hbox{if}\quad A-a< T-t,
\end{array} \right.
\end{equation}
\noindent where $\pi_1(a)=e^{-\int_0^a\mu_m(r)dr}$ and $\pi_2(a)=e^{-\int_0^a\mu_f(r)dr}.$

\noindent Suppose that the assumptions $(H_1),\;(H_2),\;(H_3)$ and $(H_4)$ are fullfiled, then we have the following result.
\begin{theorem}\label{theoIOP}
Under the assumptions of Theorem~$\ref{theo1},$ if $(a_1,a_2)\subset (b_1,b_2),$ there exists a constant $C_T>0$ independent of $p$ such that the couple $(n,l)$ solution of $(\ref{SAdjoint})$ verifies the following inequality:
\begin{equation}\label{IOP}
\displaystyle\int_0^A n^2(a,0)da+\displaystyle\int_0^A l^2(a,0)da\leq C_T\bigg(\displaystyle\int_\Xi n^2(a,t)dadt+\displaystyle\int_{\Xi'} l^2(a,t)dadt\bigg).
\end{equation}
\end{theorem}
\noindent For the proof of Theorem~$\ref{theoIOP},$ we state the following estimations of the non-local terms.

\begin{proposition}\label{Prop1}
Under the assumptions of Theorem~$\ref{theo1},$ there exists $C>0$ such that 
\begin{equation}\label{e1Prop1}
\displaystyle\int_0^{T-\eta}n^2(0,t)dt\leq C\displaystyle\int_0^T\int_{a_1}^{a_2}n^2(a,t)dadt,
\end{equation}
\noindent where $a_1<\eta<T.$

\noindent In particular, for every $\varrho>0,$ if $a_1=0$ and $n_T(a)=0$ a.e $a\in(0,\varrho);$ there is $C_{\varrho,T}>0$ such that:
\begin{equation}\label{e2Prop1}
\displaystyle\int_0^{T}n^2(0,t)dt\leq C_{\varrho}\displaystyle\int_0^T\int_{0}^{a_2}n^2(a,t)dadt.
\end{equation}
\noindent Moreover, if $(\ref{conditionrem1})$ hold, we have the inequality
\begin{equation}\label{e3Prop1}
\displaystyle\int_0^{T-\eta}l^2(0,t)dt\leq C\displaystyle\int_{\Xi'} l^2(a,t)dadt,
\end{equation}
\noindent for every $\eta$ such that $b_1<b$ and $b_1<\eta<T.$
\end{proposition}

\begin{demo}
Proposition~$\ref{Prop1}$

\noindent The state $n$ of $(\ref{SAdjoint})$ verifies 
\begin{equation}\label{e1DemoProp1}
\left\lbrace\begin{array}{l}
-\dfrac{\partial n}{\partial t}-\dfrac{\partial n}{\partial a}+\mu_m n=0\quad\hbox{in}\;(0,a_2)\times(0,T),  \\ 
\\
n(a,T)=n_T\quad\hbox{in}\;(0,a_2).
\end{array}\right. 
\end{equation}
\noindent We denote by $\tilde{n}(a,t)=n(a,t)e^{-\int_0^a\mu_m(s)ds}.$ Then, $\tilde{n}$ satisfies
\begin{equation}\label{e2DemoProp1}
\dfrac{\partial \tilde{n}}{\partial t}+\dfrac{\partial \tilde{n}}{\partial a}=0\quad\hbox{in}\;(a_1,a_2)\times(0,T).
\end{equation}
\noindent Proving the inequality
\begin{equation}\label{e3DemoProp1}
\displaystyle\int_0^{T-\eta}\tilde{n}^2(0,t)dt\leq C\int_0^T\int_{a_1}^{a_2}\tilde{n}^2(a,t)dadt
\end{equation}
\noindent lead to get inequality $(\ref{e1Prop1}).$

\noindent Indeed, we have
\begin{align*}
\displaystyle\int_0^{T-\eta}n^2(0,t)dt=\displaystyle\int_0^{T-\eta}\tilde{n}^2(0,t)dt\leq  C\int_0^T\int_{a_1}^{a_2}\tilde{n}^2(a,t)dadt & =C\int_0^T\int_{a_1}^{a_2}e^{-2\int_0^a\mu_m(s)ds}n^2(a,t)dadt\\
& \leq C\int_0^T\int_{a_1}^{a_2}n^2(a,t)dadt.
\end{align*}
\noindent We consider the following characteristics trajectory $\gamma(\lambda)=\big(T-\lambda,T+t-\lambda\big).$ If $\lambda=T$, the backward characteristic starting from $(0,t).$ If $T<a_1,$ the trajectory $\gamma(\lambda)$ never reaches the observation region $(a_1,a_2).$ So we choose $T>a_1.$

\noindent Without loss of generality, let us assume that $\eta<a_2<T.$

\begin{center}
\textbf{step 1: Estimation of} $n(0,t),\;t\in(0,T-\eta)$
\end{center}
\noindent $\bullet$ \textbf{Estimation for} $t\in(0,T-a_2).$

\vspace*{0.3cm}

\noindent We denote by:
\[w(\lambda)=\tilde{n}(T-\lambda,T+t-\lambda),\;\lambda\in(T-a_2,T).\]
\noindent Then, $\dfrac{\partial w}{\partial\lambda}=0$ for all $\lambda\in (T-a_2,T).$ In particular, $w(\lambda)$ is constant for all $\lambda\in (T-a_2,T).$ Since $\eta<a_2,$ we have
\[w(t)=\dfrac{1}{a_2-\eta}\displaystyle\int_{t-a_2}^{t-\eta}w(\lambda)d\lambda.\]
\noindent Therefore,
\[w(T)=\tilde{n}(0,t)=\dfrac{1}{a_2-\eta}\displaystyle\int_{T-a_2}^{T-\eta}\tilde{n}(T-\lambda,T+t-\lambda)d\lambda=\dfrac{1}{a_2-\eta}\displaystyle\int_{\eta}^{a_2}\tilde{n}(s,t+s)ds.\] 
\noindent Using the fact that 
\[\bigg(\displaystyle\int_{t-a_2}^{t-\eta}w(\lambda)d\lambda\bigg)^2\leq (a_2-\eta)\displaystyle\int_{t-a_2}^{t-\eta}w^2(\lambda)d\lambda\]
\noindent and integrating with respect to $t$ over $(0,T-a_2)$ we get
\[\displaystyle\int_0^{T-a_2}\tilde{n}^2(0,t)dt\leq C\displaystyle\int_{a_1}^{a_2}\int_0^{T-a_2}\tilde{n}^2(s,t+s)dtds=C\displaystyle\int_{a_1}^{a_2}\int_s^{T-a_2+s}\tilde{n}^2(s,u)duds.\]
\noindent Finally
\begin{equation}\label{e4DemoProp1}
\displaystyle\int_0^{T-a_2}\tilde{n}^2(0,t)dt\leq C\displaystyle\int_0^T\int_{a_1}^{a_2}\tilde{n}^2(a,t)dadt.
\end{equation}  

\vspace*{0.3cm}

\noindent $\bullet$ \textbf{Estimation for} $t\in(T-a_2,T-\eta)$

\vspace*{0.3cm}

\noindent We define:
\[w(\lambda)=\tilde{n}(T-\lambda,T+t-\lambda),\;\lambda\in(T-a_1,T).\] 
\noindent Then, $\dfrac{\partial w}{\partial\lambda}=0$ for all $\lambda\in (T-a_2,T).$ In particular, $w(\lambda)$ is constant for all $\lambda\in (T-a_2,T).$ Thus we have
\[w(t)=\dfrac{1}{\eta-a_1}\displaystyle\int_{t-\eta}^{t-a_1}w(\lambda)d\lambda=\dfrac{1}{\eta-a_1}\displaystyle\int_{a_1}^{\eta}\tilde{n}(s,t+s)ds.\]
\noindent Integrating with respect to $t$ over $(T-a_2,T-\eta),$ we get
\[\displaystyle\int_{T-a_2}^{T-\eta}\tilde{n}^2(0,t)dt\leq C(\eta,a_1)\displaystyle\int_{a_1}^\eta\int_{T-a_2}^{T-\eta}\tilde{n}(s,t+s)dtds.\]
\noindent Finally
\begin{equation}\label{e5DemoProp1}
\displaystyle\int_{T-a_2}^{T-\eta}\tilde{n}^2(0,t)dt\leq C(\eta,a_1)\displaystyle\int_{0}^{T}\int_{a_1}^{a_2}\tilde{n}(a,t)dadt.
\end{equation}
\noindent Combining $(\ref{e4DemoProp1})$ and $(\ref{e5DemoProp1}),$ we obtain
\[\displaystyle\int_{0}^{T-\eta}\tilde{n}^2(0,t)dt\leq C(\eta,a_1,a_2)\displaystyle\int_{0}^{T}\int_{a_1}^{a_2}\tilde{n}(a,t)dadt,\]
\noindent leading to $(\ref{e1Prop1}).$

\noindent Remark that $\displaystyle\lim_{\eta\rightarrow a_1^+}\dfrac{1}{\eta-a_1}=+\infty$ then, $\displaystyle\lim_{\eta\rightarrow a_1^+}C(\eta,a_1,a_2)=+\infty.$

\noindent Suppose now that $a_1=0.$ From the above, we have for all $\varrho>0,$ the existence of a constant depending on $\varrho$ such that:
\[\displaystyle\int_0^{T-\varrho}n^2(0,t)dt\leq C(\varrho,a_2)\displaystyle\int_0^T\int_{a_1}^{a_2}n^2(a,t)dadt.\] 
\noindent Moreover, if $n_T(a)=0$ in $(0,\varrho),$ we have $w(t)=0$ in $(T-\varrho,T).$ Then $w\equiv 0.$ So $n(0,t)=0$ in $(T-\varrho,T).$

\noindent Finally, if $n_T(a)=0$ in $(0,\varrho),$ we have the following inequality:
\[\displaystyle\int_0^{T}n^2(0,t)dt\leq C(\varrho,a_2)\displaystyle\int_0^T\int_{0}^{a_2}n^2(a,t)dadt.\]

\begin{center}
\textbf{step 2: Estimation of} $l(0,t),\;t\in(0,T-\eta)$
\end{center}
\noindent Considering $\nu=\min\{b,b_2\},$ the state $l$ of the adjoint system can be rewritten as:
\[\left\lbrace\begin{array}{l}
-\dfrac{\partial l}{\partial t}-\dfrac{\partial l}{\partial a}+\mu_f l=0\quad\hbox{in}\;(0,\nu)\times(0,T),  \\ 
\\
l(a,T)=l_T\quad\hbox{in}\;(0,\nu).
\end{array}\right.\]
\noindent We denote by $\tilde{l}(a,t)=l(a,t)e^{-\int_0^a\mu_f(s)ds}.$ Then, $\tilde{l}$ satisfies
\[\dfrac{\partial \tilde{l}}{\partial t}+\dfrac{\partial \tilde{l}}{\partial a}=0\quad\hbox{in}\;(b_1,\nu)\times(0,T).\]

\vspace*{0.3cm}

\noindent $\bullet$ \textbf{Estimation for} $t\in(0,T-\nu).$

\vspace*{0.3cm}

\noindent Defining $w$ as:
\[w(\lambda)=\tilde{l}(T-\lambda,T+t-\lambda),\;\lambda\in(T-\nu,T-b_1),\]
\noindent then, $w(\lambda)$ is constant for all $\lambda\in(T-\nu,T-b_1)$ and we have $w(t)=\dfrac{1}{\nu-\eta}\displaystyle\int_{t-\nu}^{t-\eta}w(\lambda)d\lambda.$ Likewise the step 1, we obtain
\begin{equation}\label{e6DemoProp1}
\displaystyle\int_0^{T-\nu}\tilde{l}^2(0,t)dt\leq C(\eta,\nu)\displaystyle\int_{\Xi'}\tilde{l}^2(a,t)dadt
\end{equation}

\vspace*{0.3cm}

\noindent $\bullet$ \textbf{Estimation for} $t\in(T-\nu,T-\eta).$

\vspace*{0.3cm}

\noindent We denote by:
\[w(\lambda)=\tilde{l}(T-\lambda,T+t-\lambda),\;\lambda\in(T-\eta,T-b_1).\]
\noindent As above $w(\lambda)$ is constant for all $\lambda\in(T-\eta,T-b_1)$ and we set $w(t)=\dfrac{1}{\eta-b_1}\displaystyle\int_{t-\eta}^{t-b_1}w(\lambda)d\lambda.$ Then,
\begin{equation}\label{e7DemoProp1}
\displaystyle\int_{T-\nu}^{T-\eta}\tilde{l}^2(0,t)dt\leq C(\eta,b_1)\displaystyle\int_{\Xi'}\tilde{l}^2(a,t)dadt
\end{equation}
\noindent Combining $(\ref{e6DemoProp1})$ and $(\ref{e7DemoProp1}),$ the inequality $(\ref{e3Prop1})$ follows.
\end{demo}

\begin{proposition}\label{Prop2}
Let us assume the assumptions $(H_1)-(H_3).$ For every $T>\sup\{a_1,A-a_2\}$ there exists $C_T>0$ such that the solution $(n,l)$ of the system $(\ref{EPauxiliary})$ verifies the following inequality:
\begin{equation}\label{e1Prop2}
\displaystyle\int_0^A n^2(a,0)da\leq C_T\displaystyle\int_\Xi n^2(a,t)da dt.
\end{equation}
\end{proposition}
\noindent Note that for every $T>\sup\{a_1,A-a_2\}$, there exists $a_0\in(a_1,a_2)$ such that $n(a,0)=0$ for all $a\in(a_0,A).$ This is a consequence of the following lemma.
\begin{lemma}\label{lem1}
Let us suppose that $T>\sup\{a_1,A-a_2\}.$ Then there exists $a_0\in(a_1,a_2)$ such that $T>A-a_0>A-a$ for all $a\in(a_0,A).$ Therefore, $n(a,0)=0\; \hbox{for all}\; a\in(a_0,A).$
\end{lemma}
\begin{demo}Lemma~\ref{lem1}

\noindent Suppose that $T>A-a_2,$ then there exists $\kappa>0$ (we choose $\kappa$ such that $\kappa<a_2-a_1$) $T>A-a_2+\kappa.$ So $T>A-(a_2-\kappa)$ and we denote $a_0=a_2-\kappa.$ Then, $T>A-a_0>A-a$ for all $a\in(a_0,A).$ Finally, from $(\ref{esol1})$ for all $(a,t)$ such that $T-t>A-a,$ we get $n(a,0)=0$ for all $a\in(a_0,A).$
\end{demo}
\begin{demo}Proposition~\ref{Prop2}

\noindent From the Lemma~\ref{lem1}, we have to prove the following inequality:
\begin{equation}\label{e1DemoProp2}
\displaystyle\int_0^{a_0}n^2(a,0)da\leq C_T\displaystyle\int_0^T\int_{a_1}^{a_2}n^2(a,t)dadt.
\end{equation}
\noindent We set $\tilde{n}(a,t)=e^{-\int_0^a\mu_m(s)ds}n(a,t).$ Then, from the first equation of $(\ref{SAdjoint}),\;\tilde{n}$ satisfies
\begin{equation}\label{e2DemoProp2}
\dfrac{\partial\tilde{n}}{\partial t}+\dfrac{\partial\tilde{n}}{\partial a}=0\quad\hbox{in}\;(0,A)\times(0,T).
\end{equation}
\noindent Inequality $(\ref{e2DemoProp1})$ is a consequence of the following estimation:
\[\displaystyle\int_0^{a_0}\tilde{n}^2(a,0)da\leq C\displaystyle\int_0^T\int_{a_1}^{a_2}\tilde{n}^2(a,t)dadt.\]
\noindent Indeed, we have
\begin{align*}
&\displaystyle\int_0^{a_0}n^2(a,0)da =\displaystyle\int_0^{a_0}e^{2\int_0^a\mu_m(s)ds}\tilde{n}^2(a,0)da\leq e^{2\vert\vert\mu_m\vert\vert_{L^1(0,a_0)}}\displaystyle\int_0^{a_0}\tilde{n}^2(a,0)da\leq C e^{4\vert\vert\mu_m\vert\vert_{L^1(0,a_0)}}\displaystyle\int_0^T\int_{a_1}^{a_2}n^2(a,t)dadt.
\end{align*}
\noindent We consider in this proof the characteristics $\gamma(\lambda)=(a+\lambda,\lambda)$. For $\lambda=0$ the characteristics starting from $(a,0).$ We have two cases.

\vspace*{0.3cm}

\noindent \textbf{Case 1}: $T<a_2.$

\noindent Two situations can arise:

\noindent $\bullet\;b_0=a_2-T<a_1<a_0,$ in this situation we split the interval $(0,a_0)$ as:
\[(0,a_0)=(0,b_0)\cup(b_0,a_1)\cup(a_1,a_0).\]
\noindent $\bullet\;a_1<b_0<a_0,$ in this situation we split the interval $(0,a_0)$ as:
\[(0,a_0)=(0,a_1)\cup(a_1,a_0).\]

\vspace*{0.3cm}

\noindent \textbf{Case 2}: $T\geq a_2.$

\noindent In this case, we split the interval $(0,a_0)$ as:
\[(0,a_0)=(0,a_1)\cup(a_1,a_0).\]
\noindent We make the proof for the second case. For the proof in the first case, see \cite{YSOT}.

\noindent \textbf{Upper bound on} $(0,a_1):$

\noindent For $a\in(0,a_1),$ we set $w(\lambda)=\tilde{n}(T+a-\lambda,T-\lambda),\;\lambda\in(a_1,T).$ We prove easily that $w$ is a constant, see the proof of Proposition~\ref{Prop1}. Then, we set
\[w(t)=\dfrac{1}{a_2-a_1}\displaystyle\int_{t-a_2}^{t-a_1}w(\lambda)d\lambda.\]
\noindent Integrating with respect to $a$ over $(0,a_1)$ we get
\[\displaystyle\int_0^{a_1}\tilde{n}^2(a,0)da\leq C\displaystyle\int_0^{a_1}\int_{a_1}^{a_2}\tilde{n}^2(a+\alpha,\alpha)d\alpha da.\]
\noindent Finally, we obtain
\begin{equation}\label{e3DemoProp2}
\displaystyle\int_0^{a_1}\tilde{n}^2(a,0)da\leq C\displaystyle\int_0^{T}\int_{a_1}^{a_2}\tilde{n}^2(a,t)da dt.
\end{equation} 
\noindent \textbf{Upper bound on} $(a_1,a_0):$

\noindent For $a\in(a_1,a_0),$ we set $w(\lambda)=\tilde{n}(T+a-\lambda,T-\lambda),\;\lambda\in(a_0,T)$ and 
\[w(t)=\dfrac{1}{a_2-a_0}\displaystyle\int_{t-a_2}^{t-a_0}w(\lambda)d\lambda.\]
\noindent Making as a above and integrating with respect to $a$ over $(a_1,a_0)$ it follows that
\begin{equation}\label{e4DemoProp2}
\displaystyle\int_{a_1}^{a_0}\tilde{n}^2(a,0)da\leq C\displaystyle\int_0^{T}\int_{a_1}^{a_2}\tilde{n}^2(a,t)da dt.
\end{equation} 
\noindent The inequalities $(\ref{e3DemoProp2})$ and $(\ref{e4DemoProp2})$ give the desired result.
\end{demo}
\noindent We also need the following estimate for the proof of the Theorem~\ref{theoIOP}.
\begin{proposition}\label{Prop3}
Let us assume the assumptions $(H_1)-(H_2),$ let $b_1<a_0<b$ and $T>b_1.$

\noindent Then, there exists $C_T>0$ such that the solution $l$ of $(\ref{SAdjoint})$ verifies
\begin{equation}\label{e1Prop3}
\displaystyle\int_0^{a_0}l^2(a,0)da\leq C_T\displaystyle\int_{\Xi'}l^2(a,t)da dt.
\end{equation}
\end{proposition}
\begin{demo}Proposition~\ref{Prop3}

\noindent We suppose that $\beta=0$ in $(0,b),$ the function $l$ verifies:
\[\left\lbrace\begin{array}{l}
-\dfrac{\partial l}{\partial t}-\dfrac{\partial l}{\partial a}+\mu_f l=0\quad\hbox{in}\;(0,b)\times(0,T),  \\ 
\\
l(a,T)=l_T\quad\hbox{in}\;(0,b).
\end{array}\right.\]
\noindent Proceeding as in the proof of Proposition~\ref{Prop2}, we get the desired result.
\end{demo}
\noindent For the proof of Theorem~\ref{theoIOP}, we start by the following lemma.
\begin{lemma}\label{lem2}
Suppose that $T>A-a_2+a_1$ and $a_1<b.$

\noindent Then, there exists $a_0\in(a_1,b)$ and $\kappa>0$ such that
\[T>T-(a_1+\kappa)>A-a_0>A-a\quad\hbox{for all}\;a\in(a_0,A).\]
\noindent Therefore
\[l(a,0)=\displaystyle\int_0^{A-a}\dfrac{\pi_2(a+s)}{\pi_2(a)}\beta(a+s,p(s))\big(l(0,s)+n(0,s)\big)ds\quad\hbox{for all}\;a\in(a_0,A).\]
\end{lemma}
\begin{demo}Lemma~\ref{lem2}

\noindent Notice that the solution $l$ of $(\ref{SAdjoint})$ is given by
\[l(t)=\left\lbrace\begin{array}{l}
\dfrac{\pi_2(a+T-t)}{\pi_2(a)}l_T(a+T-t)\\
+\displaystyle\int_t^T\frac{\pi_2(a+s-t)}{\pi_2(a)}\beta(a+s-t,p(s))\big((1-\gamma)n(0,s)+\gamma l(0,s)\big)ds\quad\hbox{if}\quad T-t\leq A-a, \\
\\ 
\displaystyle\int_t^{t+A-a}\frac{\pi_2(a+s-t)}{\pi_2(a)}\beta(a+s-t,p(s))\big((1-\gamma)n(0,s)+\gamma l(0,s)\big)ds \quad\hbox{if}\quad A-a< T-t.
\end{array} \right.\]
\noindent Without loss of generality, we suppose that $a_2=b.$

\noindent Suppose that $T>a_1+A-a_2\Leftrightarrow T-a_1>A-a_2.$ Then there exists $\kappa>0$ (we choose $\kappa$ such that $2\kappa<a_2-a_1$) such that $T-(a_1+\kappa)>A-(a_2-\kappa).$ We denote by $a_0=a_2-\kappa$ and as $A-a_0>A-a$ for all $a\in(a_0,A),$ then $T>T-(a_1+\kappa)>A-a_0>A-a$ for all $a\in(a_0,A).$ 

\noindent Moreover, for $(a,t)$ such that $T-t>A-a,$ we have
\[l(a,t)=\displaystyle\int_t^{t+A-a}\frac{\pi_2(a+s-t)}{\pi_2(a)}\beta(a+s-t,p(s))\big((1-\gamma)n(0,s)+\gamma l(0,s)\big)ds.\]
\noindent For $t=0$ and $a\in(a_0,A),\;T-0>A-a_0>A-a,$ one has
\[l(a,0)=\displaystyle\int_0^{A-a}\frac{\pi_2(a+s)}{\pi_2(a)}\beta(a+s,p(s))\big((1-\gamma)n(0,s)+\gamma l(0,s)\big)ds.\]
\noindent Remark that as $(a_1,a_2)\subset(b_1,b_2),$ if $a_0\in(a_1,a_2)$ then $a_0\in(b_1,b_2).$
\end{demo}
\noindent Now, we can prove the Theorem~\ref{theoIOP}.

\begin{demo}Theorem~\ref{theoIOP}

\noindent Let $a_0$ as in the Lemma~\ref{lem2}, we have:
\[\displaystyle\int_0^A l^2(a,0)da=\displaystyle\int_0^{a_0} l^2(a,0)da+\displaystyle\int_{a_0}^A l^2(a,0)da.\]
\noindent Using the results of Lemma~\ref{lem2}, the assumption of $\beta$ and the regularity of $\dfrac{\pi_2(a+s)}{\pi_2(a)},$ we can prove the existence of a constant $K_T$ independent of $p$ ($\beta\in L^\infty((0,T)\times(0,A)$) such that:
\[\displaystyle\int_{a_0}^A l^2(a,0)da\leq K_T\bigg(\displaystyle\int_0^{T-(a_1+\kappa)}n^2(0,t)dt+\displaystyle\int_0^{T-(a_1+\kappa)}l^2(0,t)dt\bigg).\]
\noindent Moreover, we have $b_1\leq a_1\leq a_1+\kappa.$ Using Proposition~\ref{Prop1}, it follows that
\[\displaystyle\int_0^{T-(a_1+\kappa)}n^2(0,t)dt\leq C_T\displaystyle\int_\Xi n^2(a,t)dadt\quad\hbox{and}\quad \displaystyle\int_0^{T-(a_1+\kappa)}l^2(0,t)\leq C_T'\displaystyle\int_{\Xi'}l^2(a,t)dadt.\]
\noindent Finally, adding the above to the results of Proposition~\ref{Prop2} and Proposition~\ref{Prop3}, we get: 
\[\displaystyle\int_0^A n^2(a,0)da+\displaystyle\int_0^A l^2(a,0)da\leq C_T\bigg(\displaystyle\int_\Xi n^2(a,t)dadt+\displaystyle\int_{\Xi'} l^2(a,t)dadt\bigg).\]
\end{demo}
\noindent For $\epsilon>0$ and $\theta>0$, we consider the functional $J_{\epsilon,\theta}$ defined by:
\begin{equation}\label{e47}
J_{\epsilon,\theta}(v_m,v_f)=\dfrac{1}{2}\displaystyle\int_\Xi v_m^2dadt+\dfrac{1}{2}\displaystyle\int_{\Xi'} v_f^2dadt+\dfrac{1}{2\epsilon}\displaystyle\int_0^A m^2(a,T)da+\dfrac{1}{2\epsilon}\displaystyle\int_0^A f^2(a,T)da,
\end{equation}
\noindent where $(m,f)$ is the solution of the following system
\begin{equation}\label{e48}
\left\lbrace\begin{array}{l}
m_t+m_a+\mu_m m=\chi_\Xi v_m\quad\hbox{in}\;Q, \\
\\ 
f_t+f_a+\mu_f f=\chi_{\Xi'} v_f\quad\hbox{in}\;Q, \\ 
\\
m(a,0)=m_0(a),\quad f(a,0)=f_0(a)\quad\hbox{in}\;Q_A, \\
\\ 
m(0,t)=(1-\gamma)\displaystyle\int_0^A\beta(a,p)f(a,t)da,\quad f(0,t)=\gamma\displaystyle\int_0^A\beta(a,p)f(a,t)da\quad\hbox{in}\;Q_T.
\end{array} \right.
\end{equation}
\begin{lemma}\label{lem3}
The functional $J_{\epsilon,\theta}$ is continuous, strictly convex and coercive. Consequently, $J_{\epsilon,\theta}$ reaches its minimum at a point $(v_{m,\epsilon},v_{f,\theta})\in L^2(\Xi)\times L^2(\Xi')$. Setting $(m_\epsilon,f_\theta)$ the associated solution of $(\ref{e48})$ and $(n_\epsilon,l_\theta)$ the solution of $(\ref{SAdjoint})$ with 
\[n_\epsilon(a,T)=-\dfrac{1}{\epsilon}m_\epsilon(a,T)\quad \hbox{and}\quad l_\theta(a,T)=-\dfrac{1}{\theta}f_\theta(a,T),\]
\noindent we have
\[\chi_\Xi v_{m,\epsilon}=\chi_\Xi n_{\epsilon}\quad\hbox{and}\quad \chi_{\Xi'} v_{f,\theta}=\chi_{\Xi'} l_{\theta}.\]
\noindent Moreover, there exit $C_i>0,\;1\leq i\leq 4,$ independent of $\epsilon$ and $\theta$ such that
\begin{align*}
\displaystyle\int_\Xi n_\epsilon^2(a,t)dadt &\leq C_1\bigg(\displaystyle\int_0^A m_0^2(a)da+\displaystyle\int_0^Af_0^2(a)da\bigg),\\
\displaystyle\int_0^A m_\epsilon^2(a,T)da &\leq\epsilon C_2\bigg(\displaystyle\int_0^A m_0^2(a)da+\displaystyle\int_0^A f_0^2(a)da\bigg),\\
\displaystyle\int_{\Xi'} l_\theta^2(a,t)dadt &\leq C_3\bigg(\displaystyle\int_0^A m_0^2(a)da+\displaystyle\int_0^Af_0^2(a)da\bigg),\\
\displaystyle\int_0^A f_\theta^2(a,T)da &\leq \theta C_4\bigg(\displaystyle\int_0^A m_0^2(a)da+\displaystyle\int_0^Af_0^2(a)da\bigg).
\end{align*}
\end{lemma}
\begin{demo}
Lemma~\ref{lem3}

\noindent It is easy to check that $J_{\epsilon,\theta}$ is coercive, continuous and strictly convex. Then, it admits a unique minimiser $(v_\epsilon,v_\theta).$ The maximum principle gives  
\begin{equation}\label{egalitecontroladjoint}
\chi_\Xi v_{m,\epsilon}=\chi_\Xi n_{\epsilon}\quad\hbox{and}\quad \chi_{\Xi'} v_{f,\theta}=\chi_{\Xi'} l_{\theta}
\end{equation}
\noindent where the couple $(n_\epsilon,l_\theta)$ is the solution of the system:
\begin{equation}\label{e49}
\left\lbrace\begin{array}{l}
-\partial_t n_\epsilon-\partial_a n_\epsilon +\mu_m n_\epsilon=0\quad\hbox{in}\;Q, \\
\\ 
-\partial_t l_\theta-\partial_a l_\theta +\mu_f l_\theta=(1-\gamma)\beta(a,p)n_\epsilon(0,t)+\gamma\beta(a,p)l_\theta(0,t)\quad\hbox{in}\;Q,\\
 \\
n_\epsilon(a,T)=-\dfrac{1}{\epsilon}m_\epsilon(a,T),\quad l_\theta(a,T)=-\dfrac{1}{\theta}f_\theta(a,T)\quad\hbox{in}\;Q_A,\\ 
 \\ 
n_\epsilon(A,t)=0,\quad l_\theta(A,t)=0 \quad\hbox{in}\;Q_T.
\end{array}\right. 
\end{equation}
\noindent Multiplying the first and the second equation of $(\ref{e49})$ by respectively $m_\epsilon$ and $f_\theta,$ integrating with respect to $Q$ and using $(\ref{egalitecontroladjoint})$ we get
\begin{align}\label{e50}
& \displaystyle\int_\Xi n_\epsilon^2(a,t)dadt+\dfrac{1}{\epsilon}\displaystyle\int_0^A m_\epsilon^2(a,T)da=-\displaystyle\int_0^A m_0(a)n_\epsilon(a,0)da-(1-\gamma)\displaystyle\int_0^T\int_0^A\beta(a,p)f_\theta(a,t)n_\epsilon(0,t)dadt
\end{align}
\noindent and 
\begin{align}\label{e51}
& \displaystyle\int_{\Xi'} l_\theta^2(a,t)dadt+\dfrac{1}{\theta}\displaystyle\int_0^A l_\theta^2(a,T)da=-\displaystyle\int_0^A f_0(a)l_\theta(a,0)da+(1-\gamma)\displaystyle\int_0^T\int_0^A\beta(a,p)f_\theta(a,t)n_\epsilon(0,t)dadt.
\end{align}
\noindent Combining $(\ref{e50})$ and $(\ref{e51}),$ we obtain
\begin{align*}
& \displaystyle\int_\Xi n_\epsilon^2(a,t)dadt+\dfrac{1}{\epsilon}\displaystyle\int_0^A m_\epsilon^2(a,T)da+\displaystyle\int_{\Xi'} l_\theta^2(a,t)dadt+\dfrac{1}{\theta}\displaystyle\int_0^A l_\theta^2(a,T)da=-\displaystyle\int_0^A m_0(a)n_\epsilon(a,0)da\\
&\hspace{11.5cm}-\displaystyle\int_0^A f_0(a)l_\theta(a,0)da.
\end{align*}
\noindent Using the inequality of Young, we have for any $\delta>0,$
\begin{align*}
& \displaystyle\int_\Xi n_\epsilon^2(a,t)dadt+\dfrac{1}{\epsilon}\displaystyle\int_0^A m_\epsilon^2(a,T)da+\displaystyle\int_{\Xi'} l_\theta^2(a,t)dadt+\dfrac{1}{\theta}\displaystyle\int_0^A l_\theta^2(a,T)da\leq\dfrac{\delta}{2}\displaystyle\int_0^A m_0^2(a)da\\
&\qquad +\dfrac{1}{2\delta}\displaystyle\int_0^A n_\epsilon^2(a,0)da +\dfrac{\delta}{2}\displaystyle\int_0^A f_0^2(a)da+\dfrac{1}{2\delta}\displaystyle\int_0^A l_\theta^2(a,0)da.
\end{align*}
\noindent Using the observability inequality $(\ref{IOP})$ and choosing $\delta=C_T$ in the previous inequality, it follows that
\begin{align*}
& \dfrac{1}{2}\displaystyle\int_\Xi n_\epsilon^2(a,t)dadt+\dfrac{1}{\epsilon}\displaystyle\int_0^A m_\epsilon^2(a,T)da+\dfrac{1}{2}\displaystyle\int_{\Xi'} l_\theta^2(a,t)dadt+\dfrac{1}{\theta}\displaystyle\int_0^A l_\theta^2(a,T)da\leq\dfrac{C_T}{2}\bigg(\displaystyle\int_0^A m_0^2(a)da\\
&\hspace{12.7cm} +\displaystyle\int_0^A f_0^2(a)da\bigg).
\end{align*}
\noindent This gives the desired result necessary to the proof of the main one.
\end{demo}
\noindent Now, we consider the system
\begin{equation}\label{e52}
\left\lbrace\begin{array}{l}
\partial_t m_\epsilon(p)+\partial_a m_\epsilon(p)+\mu_m m_\epsilon(p)=\chi_\Xi n_\epsilon\quad\hbox{in}\;Q, \\ 
\\
\partial_t f_\theta(p)+\partial_a f_\theta(p)+\mu_m f_\theta(p)=\chi_{\Xi'} l_\theta\quad\hbox{in}\;Q, \\ 
\\
m_\epsilon(p)(a,0)=m_0(a),\quad f_\theta(p)(a,0)=f_0(a)\quad\hbox{in}\;Q_A, \\
\\ 
m_\epsilon(p)(0,t)=(1-\gamma)\displaystyle\int_0^A\beta(a,p)f_\theta(p)(a,t)da,\quad f_\theta(p)(0,t)=\gamma\displaystyle\int_0^A\beta(a,p)f_\theta(p)(a,t)da \quad\hbox{in}\;Q_T,
\end{array}\right. 
\end{equation} 
\noindent where $(n_\epsilon,l_\theta)$ is the solution of $(\ref{e49})$ that minimizes the functional $J_{\epsilon,\theta}.$ We have the following result:
\begin{lemma}\label{lem4}
Under the assumptions of the Theorem~\ref{theo1}, the solutions $m_\epsilon$ and $f_\theta$ verify the following inequalities:
\begin{align}\label{e53}
& \displaystyle\int_0^A m_\epsilon^2(a,T)da+\displaystyle\int_0^T\int_0^A(1+\mu_m) m_\epsilon^2(a,t)dadt\leq C\bigg(\displaystyle\int_0^A m_0^2(a)da+\displaystyle\int_0^A f_0^2(a)da\bigg)
\end{align} 
\noindent and 
\begin{align}\label{e54}
& \displaystyle\int_0^A f_\theta^2(a,T)da+\displaystyle\int_0^T\int_0^A(1+\mu_f) f_\theta^2(a,t)dadt\leq C\bigg(\displaystyle\int_0^A m_0^2(a)da+\displaystyle\int_0^A f_0^2(a)da\bigg).
\end{align} 
\end{lemma}
\begin{demo}
Lemma~\ref{lem4}

\noindent We denote by
\[(y_\epsilon,z_\theta)=(e^{-\lambda_0 t}m_\epsilon,e^{-\lambda_0 t}f_\theta).\]
\noindent The functions $y\epsilon$ and $z_\theta$ verify
\begin{equation}\label{e55}
\partial_t y_\epsilon+\partial_a y_\epsilon+(\lambda_0+\mu_m) y_\epsilon=\chi_\Xi e^{-\lambda_0 t}n_\epsilon
\end{equation}
\noindent and
\begin{equation}\label{e56}
\partial_t z_\theta+\partial_a z_\theta+(\lambda_0+\mu_f) z_\theta=\chi_{\Xi'} e^{-\lambda_0 t}l_\theta.
\end{equation}
\noindent Multiplying the equality $(\ref{e55})$ and the equality $(\ref{e56})$ by respectively $y_\epsilon$ and $z_\theta$ and integrating with respect to $Q,$ we get 
\begin{align}\label{e57}
& \dfrac{1}{2}\displaystyle\int_0^A y_\epsilon^2(a,T)da+\dfrac{1}{2}\displaystyle\int_0^T y_\epsilon^2(A,t)dt+\displaystyle\int_0^T\int_0^A(\lambda_0+\mu_m(a)) y_\epsilon^2(a,t)dadt=\dfrac{1}{2}\displaystyle\int_0^A y_0^2(a)da\\
&+ (1-\gamma)^2\displaystyle\int_0^T\bigg(\int_0^A\beta(a,p)z_\theta da\bigg)^2dt+\displaystyle\int_0^T\int_0^A \chi_\Xi e^{-\lambda_0 t}n_\epsilon y_\epsilon dadt\nonumber
\end{align}
\noindent and 
\begin{align}\label{e58}
& \dfrac{1}{2}\displaystyle\int_0^A z_\theta^2(a,T)da+\dfrac{1}{2}\displaystyle\int_0^T z_\theta^2(A,t)dt+\displaystyle\int_0^T\int_0^A(\lambda_0+\mu_f(a)) z_\theta^2(a,t)dadt=\dfrac{1}{2}\displaystyle\int_0^A f_0^2(a)da\\
&+ \gamma^2\displaystyle\int_0^T\bigg(\int_0^A\beta(a,p)z_\theta da\bigg)^2dt+\displaystyle\int_0^T\int_0^A \chi_{\Xi'} e^{-\lambda_0 t}l_\theta z_\theta dadt.\nonumber
\end{align}
\noindent Using the Young inequality, Cauchy Schwartz inequality and the fact that $\beta$ is $L^\infty,$ we prove that:
\begin{align*}
& (1-\gamma)^2\displaystyle\int_0^T\bigg(\int_0^A\beta(a,p)z_\theta da\bigg)^2dt+\displaystyle\int_0^T\int_0^A \chi_\Xi e^{-\lambda_0 t}n_\epsilon y_\epsilon dadt \leq \vert\vert\beta\vert\vert^2_\infty\vert\vert z_\theta\vert\vert^2_{L^2(Q)}+\dfrac{1}{2}\vert\vert y_\epsilon\vert\vert^2_{L^2(Q)}+ \dfrac{1}{2}\vert\vert n_\epsilon\vert\vert^2_{L^2(\Xi)}
\end{align*}
\noindent and
\begin{align*}
& \gamma^2\displaystyle\int_0^T\bigg(\int_0^A\beta(a,p)z_\theta da\bigg)^2dt+\displaystyle\int_0^T\int_0^A \chi_{\Xi'} e^{-\lambda_0 t}l_\theta z_\theta dadt \leq \vert\vert\beta\vert\vert^2_\infty\vert\vert z_\theta\vert\vert^2_{L^2(Q)}+\dfrac{1}{2}\vert\vert z_\theta\vert\vert^2_{L^2(Q)}+ \dfrac{1}{2}\vert\vert l_\theta\vert\vert^2_{L^2(\Xi')}.
\end{align*}
\noindent Therefore, choosing $\lambda_0>(\vert\vert\beta\vert\vert^2_\infty+3/2),$ we get:
\begin{align*}
&\dfrac{1}{2}\displaystyle\int_0^A z_\theta^2(a,T)da+\displaystyle\int_0^T\int_0^A(1+\mu_f(a)) z_\theta^2(a,t)dadt\leq \dfrac{1}{2}\bigg(\vert\vert f_0\vert\vert^2_{Q_A}+\vert\vert l_\theta\vert\vert^2_{L^2(\Xi')}\bigg).
\end{align*}
\noindent Finally, applying the result of Lemma~\ref{lem3} to the above inequality, it follows that
\begin{align}\label{e60}
& \dfrac{1}{2}\displaystyle\int_0^A z_\theta^2(a,T)da+\displaystyle\int_0^T\int_0^A(1+\mu_f(a)) z_\theta^2(a,t)dadt\leq C\bigg(\displaystyle\int_0^A f_0^2(a)da+\displaystyle\int_0^A m_0^2(a)da\bigg)
\end{align}
\noindent and then the inequality $(\ref{e54})$ holds.

\noindent Likewise, we have
\begin{align*}
& \dfrac{1}{2}\displaystyle\int_0^A y_\epsilon^2(a,T)da+\displaystyle\int_0^T\int_0^A(1+\mu_m) y_\epsilon^2(a,t)dadt\leq \dfrac{1}{2}\vert\vert m_0\vert\vert^2_{Q_A}+\vert\vert\beta\vert\vert^2_\infty\vert\vert z_\theta\vert\vert^2_{L^2(Q)}+ \dfrac{1}{2}\vert\vert n_\epsilon\vert\vert^2_{L^2(\Xi)}
\end{align*}
\noindent Using the above inequality, Lemma~\ref{lem3} ant the inequality $(\ref{e60})$ we obtain
\begin{align*}
& \displaystyle\int_0^A y_\epsilon^2(a,T)da+\displaystyle\int_0^T\int_0^A(1+\mu_m) y_\epsilon^2(a,t)dadt\leq C\bigg(\displaystyle\int_0^A f_0^2(a)da+\displaystyle\int_0^A m_0^2(a)da\bigg)
\end{align*}
\noindent and then, we get the desired result.
\end{demo}
\noindent Finally, from the Lemma~\ref{lem3} and the Lemma~\ref{lem4}, if $(\epsilon,\theta)\longrightarrow(0,0)$ we get:
\[(\chi_\Xi n_\epsilon,\chi_{\Xi'} l_\theta)\rightharpoonup(\chi_\Xi v_m,\chi_{\Xi'} v_f)\;\hbox{and}\;( m_\epsilon,f_\theta)\rightharpoonup(m,f),\]
\noindent with $(m,f)$ solution of the problem $(\ref{EPauxiliary})$ and 
\[m(.,T)=f(.,T)=0\quad\hbox{a.e}\;a\in(0,A).\]

\noindent We have now the necessary ingredients for the proof of Theorem~\ref{theo1}.

\subsection{Proof of Theorem~\ref{theo1}}
In this section, we established the existence of a fixed point for the preceding auxiliary problem. Indeed, we consider that $(H_3)$ hold and we suppose to simplify that $\lambda(0)=\lambda(A)=0.$ For each $p\in L^2(0,T),$ let us denote by $\Lambda(p)\subset L^2(0,T)$ the set of all  $\displaystyle\int_0^A\lambda(a)m(p)da,$ where the couple $\big(m(p),f(p)\big)$ is the solution of the following system: 
\begin{align}\label{e321}
&\left\lbrace\begin{array}{l}
\partial_t m(p)+\partial_a m(p)+\mu_m m(p)=\chi_\Xi n\quad\hbox{in}\;Q, \\ 
\\
\partial_t f(p)+\partial_a f(p)+\mu_m f(p)=\chi_{\Xi'} l\quad\hbox{in}\;Q, \\ 
\\
m(p)(a,0)=m_0(a),\quad f(p)(a,0)=f_0(a)\quad\hbox{in}\;Q_A, \\
\\ 
m(p)(0,t)=(1-\gamma)\displaystyle\int_0^A\beta(a,p)f(p)(a,t)da,\quad f(p)(0,t)=\gamma\displaystyle\int_0^A\beta(a,p)f(p)(a,t)da \quad\hbox{in}\;Q_T,
\end{array}\right. 
\end{align} 
\noindent and $\big(n(p),l(p)\big)$ the corresponding solution of the minimizer of $J_{\epsilon,\theta}$ with $m(p)(a,T)=f(p)(a,T)=0$ for almost every $a\in(0,A)$.

\noindent We have the following result.
\begin{proposition}\label{Prop4}
Under the assumptions of the Theorem~\ref{theo1}, for any $p\in L^2(Q_T)$ the solution of problem $(\ref{e321})$  satisfies
\[\vert Y(t)\vert+\bigg\vert\bigg\vert\dfrac{d}{dt}Y\bigg\vert\bigg\vert_{L^2(0,T)}\leq C\bigg(\vert\vert m_0\vert\vert_{L^2(Q_A)}+\vert\vert f_0\vert\vert_{L^2(Q_A)}\bigg),\]
\noindent where $Y(t)=\int_0^A\lambda(a)m(p)da$ and the constant $C$ is independent of $p,\;m_0$ and $f_0.$
\end{proposition}
\begin{demo}
Proposition~\ref{Prop4}

\noindent Let $Y(t)=\displaystyle\int_0^A\lambda(a)m(p)da.$ It is easy to prove that $Y$ is solution of system
\begin{equation}\label{e1DemoProp4}
\left\lbrace\begin{array}{l}
\partial_t Y+ \displaystyle\int_0^A\mu_m(a)\lambda(a)m(p)da=R(t)\quad \hbox{in} \;Q_T,\\
\\ 
Y(0)=\displaystyle\int_0^A\lambda(a)m_0(a)da,
\end{array} \right.
\end{equation}
\noindent where
\[R(t)=\displaystyle\int_0^A\lambda'(a)m(p)da+(1-\gamma)\lambda(0)\displaystyle\int_0^A\beta(a,p)f(p)da+\displaystyle\int_0^{a_2}\lambda(a)n(p)da.\]
\noindent Using the Lemma~\ref{lem4} and the assumptions on $\beta$ and $\lambda,$ we infer that there exists $K>0$ such that
\begin{align}\label{e2DemoProp4}
&\vert\vert R\vert\vert_{L^2(Q_T)}\leq K\big(\vert\vert m_0\vert\vert_{L^2(Q_A)}+\vert\vert f_0\vert\vert_{L^2(Q_A)}\big).
\end{align}
\noindent By using $(\ref{e1DemoProp4}),$ the Young inequality and integrating on $Q_T$, we obtain
\[\displaystyle\int_0^T\vert \partial_t Y\vert^2 dt\leq 2\displaystyle\int_0^T\vert R(t)\vert^2 dt+2\displaystyle\int_0^T\bigg(\displaystyle\int_0^A\mu_m(a)\lambda(a)m_\epsilon(p)da\bigg)^2 dt.\]
Moreover, the Cauchy Schwartz inequality leads to
\begin{align*}
& \displaystyle\int_0^T\bigg(\displaystyle\int_0^A\mu_m(a)\lambda(a)m_\epsilon(p)da\bigg)^2 dt\leq \displaystyle\int_0^A\mu_m(a)\lambda(a)da \displaystyle\int_0^T\int_0^A\mu_m(a)\lambda(a)m^2_\epsilon(p)dadt. 
\end{align*} 
\noindent The inequality $(\ref{e53})$ and the fact that $\lambda\in C([0,A])$ give
\begin{align*}
& \displaystyle\int_0^T\int_0^A\mu_m(a)\lambda(a)m^2_\epsilon(p)dadt\leq K_1\bigg(\vert\vert m_0\vert\vert^2_{L^2(Q_A)}+\vert\vert f_0\vert\vert^2_{L^2(Q_A)}\bigg),
\end{align*}
\noindent where $K_1>0$ is independent of $p,\;\epsilon$ and $\theta.$ Moreover as $\lambda\mu_m\in L^1(0,A),$ and using $(\ref{e2DemoProp4}),$ it follows that
\begin{align}\label{eI1DemoProp4}
&\bigg\vert\bigg\vert\dfrac{d}{dt}Y\bigg\vert\bigg\vert_{L^2(0,T)}\leq C\bigg(\vert\vert m_0\vert\vert_{L^2(Q_A)}+\vert\vert f_0\vert\vert_{L^2(Q_A)}\bigg).
\end{align}
\noindent Now, let $\tilde{Y}=e^{-\lambda_0 t}Y.$ Then, $\tilde{Y}$ satisfies
\begin{align}\label{e3DemoProp4}
& \left\lbrace\begin{array}{l}
\partial_t \tilde{Y}+\lambda_0\tilde{Y}+ e^{-\lambda_0 t}\displaystyle\int_0^A\mu_m(a)\lambda(a)m_\epsilon(p)da=e^{-\lambda_0 t}R(t)\quad \hbox{in} \;Q_T,\\
\\ 
\tilde{Y}(0)=\displaystyle\int_0^A\lambda(a)m_0(a)da.
\end{array}\right.
\end{align}
\noindent Multiplying the first equation of $(\ref{e3DemoProp4})$ by $\tilde{Y},$ integrating on $(0,t)$ and using successively Cauchy Schwartz and Young inequalities, we deduce that
\begin{align*}
& \vert \tilde{Y}(t)\vert^2+\lambda_0\displaystyle\int_0^t\tilde{Y}^2dt\leq \vert \tilde{Y}(0)\vert^2+\displaystyle\int_0^t\tilde{Y}^2dt+\displaystyle\int_0^T\bigg(\displaystyle\int_0^A\mu_m(a)\lambda(a)m_\epsilon(p)da\bigg)^2 dt+\vert\vert R\vert\vert_{L^2(Q_T)}.
\end{align*}
\noindent Using the above calculations and choosing $\lambda_0>2,$ we get
\begin{align}\label{eI2DemoProp4}
& \vert \tilde{Y}(t)\vert^2\leq K_2\bigg(\vert\vert m_0\vert\vert_{L^2(Q_A)}+\vert\vert f_0\vert\vert_{L^2(Q_A)}\bigg).
\end{align}
\noindent The desired result comes from $(\ref{eI1DemoProp4})$ and $(\ref{eI2DemoProp4}).$
\end{demo}
\noindent It is obvious that $\Lambda(p)$ is convex, and let 
\[W(0,T)=\bigg\{Y\in L^\infty(0,T),\;\vert\vert Y\vert\vert_{L^\infty(0,T)}\leq N;\;\bigg\vert\bigg\vert\dfrac{dY}{dt}\bigg\vert\bigg\vert_{L^2(0,T)}\leq N\bigg\}, \]
\noindent with $N=C\bigg(\vert\vert m_0\vert\vert_{L^2(Q_A)}+\vert\vert f_0\vert\vert_{L^2(Q_A)}\bigg).$

\noindent We have $W(0,T)\subset W^{1,1}(0,T).$ Moreover the injection of $W^{1,1}(0,T)$ into $L^2(0,T)$ is compact, see \cite{HBrezis} Page $129.$ So $W(0,T)$ is relatively compact in $L^2(0,T)$. From Proposition~\ref{Prop4} we have $\Lambda\big(W(0,T)\big)\subset W(0,T),$ and we see that $\Lambda\big(W(0,T)\big)$ is a relatively compact subset of $L^2(0,T).$ Let us now prove that $\Lambda$ is upper-semicontinuous. This is equivalent to prove that for any closed subset $K$ of $L^2(0,T)$, $\Lambda^{-1}(K)$ is closed in $L^2(0,T).$ Let $(p_k)\in\Lambda^{-1}(K)$ such that $p_k$ converges towards $p$ in $L^2(0,T).$ Then, $p_k$ is bounded and for all $k$ there exists $P_k\in K$  such that $P_k\in\Lambda(p_k).$ Therefore, from the definition of $\Lambda,$ there exists $(m_k,f_k)\in \big(L^2((0,T)\times(0,A))\big)^2$ associated to $(n_k,l_k)\in L^2(\Xi)\times L^2(\Xi')$ solution of $(\ref{e321})$ such that $P_k=\displaystyle\int_0^A\lambda(a)m_k(p_k)da$ and  satisfying the inequalities of the Lemma~\ref{lem3} and Lemma~\ref{lem4}. Consequently $(m_k,f_k)$ and $(n_k,l_k)$ are bounded respectively in $\big(L^2((0,T)\times(0,A))\big)^2$ and in $L^2(\Xi)\times L^2(\Xi')$. Thus, there exists a subsequences still denoted by $(m_k,f_k)$ and $(n_k,l_k)$ that converge weakly to $(m,f)$ in $\big(L^2((0,T)\times(0,A))\big)^2$  and $(n,l)$ in $L^2(\Xi)\times L^2(\Xi')$ respectively. Using hypothesis $(H_3),$ it follows that $\displaystyle\int_0^A\lambda(a)m_k(p_k)da$ converges strongly to $\displaystyle\int_0^A\lambda(a)m(p)da$ in $L^2(0,T).$

\noindent Now, by standard device we see that  $(m,f)$ associated to $(n,l)$ are solution of $(\ref{e321})$ and  satisfy the inequalities of the Lemma~\ref{lem3} and Lemma~\ref{lem4}. This implies  that $P\in\Lambda(p).$

\noindent On the other hand, thanks to the Proposition~\ref{Prop4}, one can extract a
subsequence also denoted by $P_k$ that converges strongly towards the function $P$ in $L^2(0,T)$.
Since $K$ is closed we deduce that $P\in K$. Finally, we deduce that $p\in\Lambda^{-1}(P)$.

\noindent Applying the Kakutani fixed point theorem \cite{KD} in the space $L^2(0,T)$ to the mapping $\Lambda,$ we infer that there is at least one $Y\in W(0,T)$ such that $Y\in\Lambda(Y).$ This completes the null controllability proof of the model $(\ref{EP}).$

\subsection{Proof of Theorem~\ref{theo2}}

\subsubsection{Proof of Theorem~\ref{theo2}-(1)}
In this section, we always consider the following system:
\begin{equation}\label{e3311}
\left\lbrace\begin{array}{l}
m_t+m_a+\mu_m m=\chi_\Theta v_m\quad\hbox{in}\;Q, \\ 
\\
f_t+f_a+\mu_f f=0\quad\hbox{in}\;Q, \\
\\ 
m(a,0)=m_0,\quad f(a,0)=f_0 \quad\hbox{in}\;Q_A,\\
\\ 
m(0,t)=(1-\gamma)\int_0^A\beta(a,p)fda,\quad f(0,t)=\gamma\int_0^A\beta(a,p)fda \quad\hbox{in}\;Q_T,
\end{array}\right. 
\end{equation}
\noindent for every $p$ in $L^2(Q_T).$ Under the assumptions of Theorem~\ref{theo2}, the controllability problem that is to find $v_m\in L^2(\Theta)$ such that $(m,f)$ solution of the system $(\ref{e3311})$ verifies
\[m(.,T)=0\quad a\in(\varrho,A)\]
\noindent is equivalent to the following observability inequality.
\begin{proposition}\label{Prop3311}
Let us assume the assumptions $(H_1)-(H_2)-(H_3),$ for every $T>A-a_2$ and for any $\varrho>0,$ if $h_T(a)=0\;a.e$ in $(0,\varrho),$ there exists $C_{\varrho,T}>0$ such that the following inequality
\begin{equation}\label{e3312}
\displaystyle\int_0^A h^2(a,0)da+\displaystyle\int_0^A g^2(a,0)da\leq C_{\varrho,T}\displaystyle\int_\Theta h^2(a,t)dadt
\end{equation}
\noindent holds, where $(h,g)$ is the solution of 
\begin{equation}\label{e3313}
\left\lbrace\begin{array}{l}
-h_t-h_a+\mu_m h=0\quad\hbox{in}\;Q, \\ 
\\
-g_t-g_a+\mu_f g=(1-\gamma)\beta(a,p)g(0,t)+\gamma\beta(a,p)h(0,t)\quad\hbox{in}\;Q, \\ 
\\
h(a,T)=h_T,\quad g(a,T)=0\quad\hbox{in}\;Q_A, \\
\\ 
h(A,t)=0,\quad g(A,t)=0\quad\hbox{in}\;Q_T.
\end{array} \right.
\end{equation}
\end{proposition}
\noindent For the proof of the Proposition~\ref{Prop3311}, we state the following estimate.
\begin{proposition}\label{Prop3312}
Under the assumptions $(H_1)$ and $(H_2),$ there exists a constant $C>0$ such that the solution $(h,g)$ of the system $(\ref{e3313})$ verifies
\begin{equation}\label{e3314}
\displaystyle\int_0^A g^2(a,0)da+\displaystyle\int_0^T\int_0^A(1+\mu_f) g^2(a,t)dadt\leq C\displaystyle\int_0^T h^2(0,t)dt.
\end{equation}
\noindent Moreover, we deduce for $h_T=0\;a.e$ in $(0,\varrho)$ that there exists a constant $C_{\varrho,T}>0$ such that
\begin{equation}\label{e3315}
\displaystyle\int_0^A g^2(a,0)da\leq C_{\varrho,T}\displaystyle\int_0^T h^2(0,t)dt\leq C_{\varrho,T}\displaystyle\int_0^T\int_0^{a_2} h^2(a,t)dadt.
\end{equation}
\end{proposition}

\begin{demo}
Proposition~\ref{Prop3312}

\noindent Setting $y=e^{\lambda_0 t}g,$ the function $y$ verifies
\begin{equation}\label{e1DemoProp3312}
-\partial_t y-\partial_a y+(\lambda_0+\mu_f)y=\gamma\beta(a,p)e^{\lambda_0 t}h(0,t)+(1-\gamma)\beta(a,p)y(0,t).
\end{equation}
\noindent Multiplying the equality $(\ref{e1DemoProp3312})$ by $y$ and integrating on $Q$, we obtain:
\begin{align}\label{e2DemoProp3312}
& \dfrac{1}{2}\displaystyle\int_0^A y^2(a,0)da+\dfrac{1}{2}\displaystyle\int_0^T y^2(0,t)dt+\displaystyle\int_0^T\int_0^A (\lambda_0+\mu_f)y^2(a,t)dadt=\displaystyle\int_0^T\int_0^A\bigg(\gamma\beta(a,p)e^{\lambda_0 t}h(0,t)+(1-\gamma)\beta(a,p)y(0,t)\bigg)y dadt.
\end{align}
\noindent Using Young inequality and the condition on $\beta,$ we get
\begin{align*}
&\displaystyle\int_0^T\int_0^A\bigg(\gamma\beta(a,p)e^{\lambda_0 t}h(0,t)+(1-\gamma)\beta(a,p)y(0,t)\bigg)y dadt\leq \dfrac{e^{2\lambda_0 T}\vert\vert\beta\vert\vert^2_\infty}{2\delta}\vert\vert h(0,.)\vert\vert^2_{L^2(Q_T)}+\dfrac{\delta}{2}\vert\vert y\vert\vert^2_{L^2(Q)}+\dfrac{\vert\vert\beta\vert\vert^2_\infty}{2\delta}\vert\vert y(0,.)\vert\vert^2_{L^2(Q_T)}\\
&\hspace{9.5cm} +\dfrac{\delta}{2}\vert\vert y\vert\vert^2_{L^2(Q)}.
\end{align*}
\noindent Choosing $\delta=\vert\vert\beta\vert\vert^2_{\infty},$ we obtain
\begin{align*}
& \dfrac{1}{2}\displaystyle\int_0^A y^2(a,0)da+\displaystyle\int_0^T\int_0^A (\lambda_0+\mu_f)y^2(a,t)dadt\leq \dfrac{e^{2\lambda_0 T}}{2}\vert\vert h(0,.)\vert\vert^2_{L^2(Q_T)}+\vert\vert\beta\vert\vert^2_{\infty}\vert\vert y\vert\vert^2_{L^2(Q)}.
\end{align*}
\noindent Finally, choosing $\lambda_0>\vert\vert\beta\vert\vert^2_{\infty}+1$ it follows that
\begin{align*}
& \dfrac{1}{2}\displaystyle\int_0^A y^2(a,0)da+\displaystyle\int_0^T\int_0^A (1+\mu_f)y^2(a,t)dadt\leq \dfrac{e^{2\lambda_0 T}}{2}\displaystyle\int_0^A h^2(0,t)dt.
\end{align*}
\noindent So,
\begin{align*}
& \displaystyle\int_0^A y^2(a,0)da\leq e^{2\lambda_0 T}\displaystyle\int_0^A h^2(0,t)dt.
\end{align*}
\noindent Finally, we get
\begin{align*}
& \displaystyle\int_0^A g^2(a,0)da\leq e^{2\lambda_0 T}\displaystyle\int_0^A h^2(0,t)dt.
\end{align*}
\noindent Combining the above inequality and the inequality $(\ref{e2Prop1})$ of the Proposition~\ref{Prop1}, for $h_T=0$ a.e in $(0,\varrho)$ we get
\begin{align}\label{e3DemoProp3312}
&\displaystyle\int_0^A g^2(a,0)da\leq C_{\varrho,T}\displaystyle\int_0^T\int_0^{a_2} h^2(a,t)dadt.
\end{align}
\end{demo}

\begin{demo}
Proposition~\ref{Prop3311}

\noindent We use the results of Proposition~\ref{Prop2} and Proposition~\ref{Prop3312}. Indeed, combining $(\ref{e1Prop2})$ and $(\ref{e3DemoProp3312})$ the desired result is obtained.
\end{demo}
\noindent Now, let $\epsilon>0$ and $\varrho>0.$ We consider the functional $J_\epsilon$ defined by
\begin{equation}\label{Jepsilon1}
J_\epsilon(v_m)=\dfrac{1}{2}\displaystyle\int_0^T\int_{a_1}^{a_2}v_m^2(a,t)dadt+\dfrac{1}{2\epsilon}\displaystyle\int_{\varrho}^A m^2(a,T)da,
\end{equation}
\noindent where $(m,f)$ is the solution of the following system
\begin{equation}\label{e3316}
\left\lbrace\begin{array}{ll}
m_t+m_a+\mu_m m  =\chi_\Xi v_m &\hbox{in}\;Q, \\
\\ 
f_t+f_a+\mu_f f  =0 &\hbox{in}\;Q, \\
\\ 
m(a,0)=m_0(a)\quad f(a,0)=f_0(a) & \hbox{in}\;Q_A, \\
\\ 
m(0,t)=(1-\gamma)\displaystyle\int_0^A\beta(a,p)f(a,t)da & \hbox{in}\;Q_T, \\
\\ 
f(0,t)=\gamma\displaystyle\int_0^A\beta(a,p)f(a,t)da & \hbox{in}\;Q_T.
\end{array}\right. 
\end{equation}
\noindent We have the following lemma.
\begin{lemma}\label{lem3311}
The functional $J_\epsilon$ is continuous, strictly convex and coercive. Consequently, $J_\epsilon$ reaches its minimum at one has $v_{m,\epsilon}=\chi_\Theta h_\epsilon$ and there exists a positive constants $C_1,\;C_2$ independent of $\epsilon$ such that
\begin{align*}
& \displaystyle\int_0^T\int_0^{a_2}h_\epsilon^2(a,t)dadt\leq C_1\bigg(\displaystyle\int_0^A m_0^2(a)da+\displaystyle\int_0^A f_0^2(a)da\bigg)
\end{align*} 
\noindent and
\begin{align*}
&  \displaystyle\int_{\varrho}^A m_\epsilon^2(a,T)da\leq \epsilon C_2\bigg(\displaystyle\int_0^A m_0^2(a)da+\displaystyle\int_0^A f_0^2(a)da\bigg).
\end{align*}
\end{lemma}
\begin{demo}
Lemma~\ref{lem3311}

\noindent The proof is similar to that of Lemma~\ref{lem3}.
\end{demo}
\noindent By making $\epsilon$ tending towards zero, we thus obtain that $\chi_\Theta h_\epsilon\rightharpoonup \chi_\Theta v_m$ and $(m\epsilon,f_\epsilon)\rightharpoonup(m,f),$ where $(m,f)$ is the solution of the system $(\ref{e3316})$ that verifies
\[m(.,T)=0\;\hbox{a.e in}\;(\varrho,A).\]
\noindent Finally, a similar function $\Lambda$ is defined and a similar procedure is followed to get the null controllability for the nonlinear problem.

\subsubsection{Proof of Theorem~\ref{theo2}-(2)}
Let $p\in L^2(Q_T),$ under the assumptions of Theorem~\ref{theo2}, the following controllability problem find $v_f\in L^2(\Theta)$ such that the solution of the system
\begin{equation}\label{e3321}
\left\lbrace\begin{array}{ll}
f_t+f_a+\mu_f f  =\chi{\Xi'}v_f &\hbox{in}\;Q, \\
\\ 
f(a,0)=f_0(a) & \hbox{in}\;Q_A, \\
\\ 
f(0,t)=\gamma\displaystyle\int_0^A\beta(a,p)f(a,t)da & \hbox{in}\;Q_T
\end{array}\right. 
\end{equation}

\noindent verifies
\[f(.,T)=0\;\hbox{a.e in}\;(0,A).\]
\noindent is equivalent to the following observability inequality.
\begin{proposition}\label{Prop3321}
Let us assume true the assumptions $(H_1)-(H_2)-(H_3).$ For any $T>a_1+A-a_2$ there exists $C_T>0$ such that
\begin{equation}\label{e1Prop3321}
\displaystyle\int_0^A g^2(a,0)da\leq C_T \displaystyle\int_{\Xi'}g^2(a,t)dadt,
\end{equation}
\noindent where $g$ is solution of the system
\begin{equation}\label{e2Prop3321}
\left\lbrace\begin{array}{l}
-g_t-g_a+\mu_f g=\gamma\beta(a,p)g(0,t)\quad\hbox{in}\;Q, \\ 
\\
g(a,T)=g_T\quad\hbox{in}\;Q_A, \\
\\ 
g(A,t)=0\quad\hbox{in}\;Q_T.
\end{array} \right.
\end{equation}
\end{proposition}

\begin{demo}Proposition~\ref{Prop3321}

\noindent Using the inequality $(\ref{e3Prop1})$ of Proposition~\ref{Prop1}, the result of Proposition~\ref{Prop3} and the representation of the solution of the system $(\ref{e2Prop3321}),$ we get the desired result.
\end{demo}
\noindent To conclude, a similar function $\Lambda$ is defined and a similar procedure is followed to get the null controllability for the nonlinear system. We omit all details because the extension is straightforward.

\section*{Conclusion}
In this paper, we have proved the null controllability of a nonlinear age and two-sex structured population dynamics
model. We considered two controllability issues.

\noindent The first problem is related to the total extinction, which means that, we have estimated a time T to bring the male and female subpopulation density to zero.

\noindent The second deals with the null controllability of the density of male or female individuals.
In this case, the control is made to act either on the males or on the females. In the case where
control acts on males, we show that the density of male individuals can be reduced to zero over part of the age range. In the event that the control acts on the females, the density of female individuals can be reduced to zero over the entire area
ages.

\noindent In the desire to consolidate even more the result obtained in this article, we will think in our future work of the case of controllability with positivity constraint and also to establish algorithms allowing to calculate the control from its characterization.
These outstanding questions are therefore in this order:
\begin{enumerate}
\item[$\bullet$] \textbf{Controllability with positivity constraints}: We will be interested in controllability with positivity constraint.
\item[$\bullet$] \textbf{Numerical implementation}: For a given fertility rate $\beta$ and $\lambda$, the mortality rate $\mu_m$ and $\mu_f$, the initial conditions $m_0,\;f_0$ and a
positive parameters $\epsilon,\;\theta>0$, how to determine a numerical algorithm allowing to determine the $(\epsilon,\theta)-$approximate null
control functions $v_m$ and $v_f$?
\end{enumerate}

\end{document}